\documentclass[twoside,11pt]{article}

\usepackage{jmlr2e}
\usepackage{bm, amsmath, amsfonts, amssymb}
\usepackage{multirow}
\usepackage[figuresright]{rotating}
\usepackage{amsmath}
\usepackage{cases}
\usepackage{color}
\usepackage{rotating}
\usepackage{setspace}
\usepackage{graphics, graphicx}
\usepackage{subfigure}
\usepackage{float}
\usepackage{epsfig}
\usepackage{lscape}
\usepackage{multirow}
\usepackage{mathrsfs,amssymb}
\usepackage{lscape}
\usepackage{threeparttable}
\usepackage{bm, amsmath, amsfonts, amssymb}
\usepackage{natbib}
\usepackage{hyperref}

\usepackage{array,booktabs}

\newcommand{\beq}{\begin{equation}}
\newcommand{\eeq}{\end{equation}}
\newcommand{\beas}{\begin{align*}}
\newcommand{\eeas}{\end{align*}}
\newcommand{\bea}{\begin{align}}
\newcommand{\eea}{\end{align}}
\newcommand{\bet}{\begin{theorem}}
	\newcommand{\eet}{\end{theorem}}
\newcommand{\bel}{\begin{lemma}}
	\newcommand{\eel}{\end{lemma}}
\newcommand{\bep}{\begin{proposition}}
	\newcommand{\eep}{\end{proposition}}
\newcommand{\R}{\mathbb{R}}
\newcommand{\E}{\mathbb{E}}

\newcommand{\bu}{\bold{u}}

\newcommand{\bv}{\bold{v}}
\newcommand{\bg}{\bold{g}}

\newcommand{\bh}{\bold{h}}
\newcommand{\bU}{\bold{U}}
\newcommand{\bW}{\bold{W}}
\newcommand{\bV}{\bold{V}}
\newcommand{\bG}{\bold{G}}
\newcommand{\bQ}{\bold{Q}}
\newcommand{\bH}{\bold{H}}
\newcommand{\bZ}{\bold{Z}}
\newcommand{\bY}{\bold{Y}}
\newcommand{\bX}{\bold{X}}
\newcommand{\bD}{\bold{D}}
\newcommand{\bO}{\bold{O}}
\newcommand{\bGam}{\bold{\Gamma}}
\newcommand{\bSig}{\bold{\Sigma}}

\newcommand{\bmu}{\boldsymbol{\mu}}
\newcommand{\btheta}{\boldsymbol{\theta}}

\newcommand{\cP}{\mathcal{P}}

\newcommand{\CC}{\mathcal{C}}
\newcommand{\bS}{\mathbb{S}}
\newcommand{\bB}{\mathbb{B}}
\newcommand{\supp}{\text{supp}}

\def\PP{{\mathbb P}}

\newcommand{\argmin}{\mathop{\rm arg\min}}
\newcommand{\argmax}{\mathop{\rm arg\max}}

\jmlrheading{1}{2000}{1-43}{00/00}{00/00}{}{}

\ShortHeadings{Optimal Structured Principal Subspace Estimation}{Cai, Li and Ma}
\firstpageno{1}

\begin{document}

\title{Optimal Structured Principal Subspace Estimation: Metric Entropy and Minimax Rates}

\author{\name Tony Cai \email tcai@wharton.upenn.edu \\
       \addr Department of Statistics\\
       University of Pennsylvania\\
      Philadelphia, PA 19104, USA  
       \AND       
	   \name Hongzhe Li \email hongzhe@pennmedicine.upenn.edu \\
	   \addr Department of Biostatistics, Epidemiology and Informatics\\
	   University of Pennsylvania\\
	   Philadelphia, PA 19104, USA
	      \AND   
   \name Rong Ma \email rongm@upenn.edu \\
   \addr Department of Biostatistics, Epidemiology and Informatics\\
   University of Pennsylvania\\
   Philadelphia, PA 19104, USA}

 \editor{}

\maketitle

\begin{abstract}
Driven by a wide range of applications, several principal subspace estimation problems have been studied individually under different structural constraints. This paper presents a unified framework for the statistical analysis of a general structured principal subspace estimation problem which includes as special cases sparse PCA/SVD, non-negative PCA/SVD, subspace constrained PCA/SVD, and spectral clustering.  General minimax lower and upper bounds are established to characterize the interplay between the information-geometric complexity of the structural set for the principal subspaces, the signal-to-noise ratio (SNR), and the dimensionality. The results yield interesting phase transition phenomena concerning the rates of convergence as a function of the SNRs and the fundamental limit for consistent estimation. Applying the general results to the specific settings yields the minimax rates of convergence for those problems, including the previous unknown optimal rates for sparse SVD, non-negative PCA/SVD and subspace constrained PCA/SVD. 
\end{abstract}


\section{Introduction}

Spectral methods such as the principal component analysis (PCA) and singular value decomposition (SVD) are a ubiquitous technique in modern data analysis with a wide range of applications in many fields including statistics, machine learning, applied mathematics, and engineering. As a fundamental tool for dimension reduction, the spectral methods aim to extract the low-dimensional structures embedded in the high-dimensional data. In many of these modern applications, the complexity of the datasets and the need of incorporating the existing knowledge from the subject areas require the data analysts to take into account the prior structural information on the statistical objects of interest in their analysis. In particular, many interesting problems in high-dimensional data analysis can be formulated as a structured principal subspace estimation problem where one has the prior knowledge that the underlying principal subspace satisfies certain structural conditions (see Section \ref{related.works} for a list of related problems). 

The present paper aims to provide a unified treatment of the structured principal subspace estimation problems that have attracted much recent interest in both theory and practice.

\subsection{Problem Setup} \label{sec1.1}
To fix ideas, we consider two generic models that have been extensively studied in the literature, namely, the \emph{matrix denoising model} and the \emph{spiked Wishart model} (see, for example, \cite{johnstone2001distribution,baik2006eigenvalues,paul2007asymptotics,bai2008central,cai2013sparse,donoho2014minimax,wang2017asymptotics,choi2017selecting,donoho2018optimal,perry2018optimality,bao2018singular}, among many others). 

\begin{definition}[Matrix Denoising Model]
	Let $\bY\in \R^{p_1\times p_2}$ be the observed data matrix generated from the model $\bY= \bU\bGam\bV^\top+\bZ$ where $\bZ\in \R^{p_1\times p_2}$ has i.i.d. entries from $N(0,\sigma^2)$, $\bGam\in \R^{r\times r}$ is a diagonal matrix with ordered diagonal entries $\lambda_1\ge \lambda_2\ge...\ge \lambda_r> 0$ for $1\le r\le \min\{p_1,p_2\}$, $\bU\in O(p_1,r)$, and $\bV\in O(p_2,r)$ with $O(p,r)=\{\bW\in \R^{p\times r}:\bW^\top\bW={\bf I}_r\}$ being the set of all $p\times r$ orthonormal matrices.
\end{definition}

\begin{definition}[Spiked Wishart Model]
	Let $\bY\in \R^{n\times p}$ be the observed data matrix whose  rows $Y_i\in\R^{p}$, $i=1, \ldots, n$, are independently generated from $N(\bmu,  \bU\bGam \bU^\top+\sigma^2{\bf I}_{p})$ where $\bU\in O(p,r)$ with $1\le r\le p$, and $\bGam\in \R^{r\times r}$ is diagonal with ordered diagonal entries $\lambda_1\ge...\ge\lambda_r> 0$. Equivalently,  $Y_i$ can be viewed as $Y_i=X_i+\epsilon_i$ where $X_i \sim N(\bmu,  \bU\bGam\bU^\top)$, $\epsilon_i \sim N(0,\sigma^2{\bf I}_{p})$, and  $X_1, \ldots, X_n$ and $\epsilon_1, \ldots, \epsilon_n$ are independent. 
\end{definition}

In the past decades, these two models have attracted substantial practical and theoretical interest and have been studied in different contexts in statistics, probability, and machine learning.
This paper addresses  the problem of optimal estimation of the principal (eigen/singular) subspaces spanned by the orthonormal columns of $\bU$  (denoted as $\text{span}(\bU)$),  based on the data matrix $\bY$ and the prior structural knowledge on $\bU$. Specifically,  we aim to uncover the deep connections between the statistical limit of the estimation problem  as measured by the minimax risk and the geometric complexity of the parameter spaces as characterized by functions of certain entropy measures. 

Since the principal subspaces can be uniquely identified with their associated projection matrices, estimating $\text{span}(\bU)$ is equivalent to estimating $\bU\bU^\top$. A commonly used metric for gauging the distance between two linear subspaces $\text{span}(\bU_1)$ and $\text{span}(\bU_2)$ is
\[
d(\bU_1,\bU_2)=\|\bU_1\bU_1^\top-\bU_2\bU_2^\top\|_F.
\]
In this paper, we use $d(\cdot, \cdot)$ as the loss function and measure the performance of an estimator $\widehat{\bU}$ of $\bU$ by the risk
\[
\mathcal{R}(\widehat{\bU},\bU)=\E d(\widehat{\bU},\bU).
\]

\subsection{Related Works} \label{related.works}

The problem considered in this paper can be viewed as a generalization and unification of many interesting problems in high-dimensional statistics and machine learning. We first present a few examples to demonstrate the richness of the structured principal subspace estimation problem and its connections to the existing literature.

\begin{enumerate}
	
	\item \emph{Sparse PCA/SVD.}
	The goal of sparse PCA/SVD is to recover $\text{span}(\bU)$ under the assumption that columns of $\bU$ are sparse.
	Sparse PCA has been extensively studied in the past two decades under the spiked Wishart model  (see, for example, \cite{d2005direct,zou2006sparse,shen2008sparse,witten2009penalized,yang2011sparse,vu2012minimax,cai2013sparse,ma2013sparse,birnbaum2013minimax,Cai2015spiked}, among many others). In particular, the exact minimax rates of convergence under the loss $d(\cdot,\cdot)$ was established by \cite{cai2013sparse} in the general rank-$r$ setting. In contrast, theoretical analysis for the sparse SVD is relatively scarce, and the minimax rate of convergence remains unknown.
	
		\item \emph{Non-negative PCA/SVD.}
	Non-negative PCA/SVD aims to estimate $\text{span}(\bU)$ under the assumption that entries of $\bU$ are non-negative.
	This problem has been studied by \cite{deshpande2014cone} and \cite{montanari2015non} under the rank-one matrix denoising model ($r$=1), where the statistical limit and certain sharp asymptotics were carefully established. However, it is still unclear what are the minimax rates of convergence for estimating $\text{span}(\bU)$ under either rank-one or general rank-$r$ settings under either the spiked Wishart model or matrix denoising model.
	
	\item \emph{Subspace Constrained PCA/SVD.}
	The subspace constrained PCA/SVD  assumes the columns of $\bU$ are in some low-dimensional linear subspaces of $\R^p$. In other words, $\bU\in\CC_A(p,k)=\{ \bU\in O(p,r): A\bU_{.j}=0\text{ for all $1\le j\le r$} \}$ for some rank $(p-k)$ matrix $A\in \R^{p\times (p-k)}$ where $r< k< p$. Estimating the principal subspaces under various linear subspace constraints has been considered in many applications such as network clustering \citep{wang2010flexible,kawale2013constrained,kleindessner2019guarantees}. However, the minimax rates of convergence for subspace constrained PCA/SVD remain unknown.
	
	\item \emph{Spectral Clustering.}
	Suppose we observe $Y_i\sim N(\btheta_i,\sigma^2{\bf I}_p)$ independently, where $\btheta_i\in\{ \btheta,-\btheta\}\subset \R^p$ for $i=1,...,n$. Let $\bY\in \R^{n\times p}$ such that $Y_i$ is the $i$-th row of $\bY$. We have $\bY=\bh \btheta^\top+\bZ$ where $\bh\in \{\pm 1\}^n$ and $\bZ$ has i.i.d. entries from $N(0,\sigma^2)$. Spectral clustering of $\{Y_i\}_{1\le i\le n}$ aims to recover the class labels in $\bh$. Equivalently, spectral clustering can be treated as estimating the leading left singular vector $\bu=\bh/\|\bh\|_2$ in the matrix denoising model with $\bu\in\CC^n_{\pm}=\{\bu\in \R^n:\|\bu\|_2=1, u_i\in\{\pm n^{-1/2}\}\}$. See  \cite{azizyan2013minimax,jin2016influential,lu2016statistical,jin2017phase,cai2018rate,giraud2018partial,ndaoud2018sharp,loffler2019optimality} and references therein for recent theoretical results.
\end{enumerate}

In addition to the aforementioned problems, there are many other interesting problems that share the same generic form as the structured principal subspace estimation problem. For example, motivated by applications in the statistical analysis of metagenomics data, \cite{ma2019optimalb,ma2019optimala} considered an approximately rank-one matrix denoising model where the leading singular vector satisfies the monotonicity constraint. In a special case of matrix denoising model, namely, the Gaussian Wigner model  $\bY=\lambda\bu\bu^\top+\bZ\in\R^{n\times n}$, where $\bZ$ has i.i.d. entries (up to symmetry)  drawn from a Gaussian distribution, the { Gaussian $\mathbb{Z}/2$ synchronization} problem \citep{javanmard2016phase,perry2018optimality} aims to recover the leading singular vector $\bu$ where $\bu\in \{\bu\in \R^n:\|\bu\|_2=1, u_i\in\{\pm n^{-1/2}\}\}$. These important applications  provide motivations for a unified framework to study the fundamental difficulty and optimality of these estimation problems.

On the other hand, investigations of metric entropy as a measure of statistical complexity has been one of the central topics in theoretical statistics, ranging from nonparametric function estimation \citep{yatracos1988lower,haussler1997mutual,yang1999information,yang1999minimax,wu2016minimax}, high-dimensional statistical inference \citep{raskutti2011minimax,verzelen2012minimax,vu2012minimax,cai2013sparse,ma2013sparse} to statistical learning theory \citep{haussler1997metric,lugosi1999adaptive,bousquet2002some,bartlett2002rademacher,koltchinskii2006local,lecue2009aggregation,cai2016geometric,rakhlin2017empirical}. Among them, interesting connections between the complexity of the parameter space and the fundamental difficulty of the statistical problem as quantified by certain minimax risk have been carefully established. In this sense, the current work stands as a step along this direction in the context of principal subspace estimation under some general random matrix models.

\subsection{Main Contribution}

The main contribution of this paper is three-fold. Firstly, a unified framework is introduced for the study of structured principal subspace estimation problems under both the matrix denoising model and the spiked Wishart model. Novel generic minimax lower bounds and risk upper bounds are established to characterize explicitly the interplay between the information-geometric complexity of the structural set for the principal subspaces, the signal-to-noise ratio (SNR), and the dimensionality of the parameter spaces. The results yield interesting phase transition phenomena concerning the rates of convergence as functions of the SNRs and the fundamental limit for consistent estimation. The general lower and upper bounds reduce determination of the minimax optimal rates for many interesting problems to mere calculations of certain information-geometric quantities. 
Secondly, to obtain the general risk upper bounds, new technical tools are developed for the analysis of the proposed estimators in their general forms. In addition, the minimax lower bounds rely on careful constructions of multiple composite hypotheses about the structured parameter spaces, and non-trivial calculations of the Kullback-Leibler (KL) divergence between certain mixture probability measures, which can be of independent interest. 
Thirdly, by directly applying our general results to the specific problems discussed in Section \ref{related.works},  we establish the minimax optimal rates for those problems. Among them,  the minimax rates for  sparse SVD, non-negative PCA/SVD and subspace constrained PCA/SVD, are to our knowledge previously unknown. 

\subsection{Organization and Notation}

The rest of the paper is organized as follows. After introducing the notation at the end of this section,  we characterize in Section \ref{mdm.minimax.sec} a minimax lower bound under the matrix denoising model using local metric entropy measures. A general estimator is introduced in Section \ref{mle.sec} and its risk upper bound is obtained via certain global metric entropy measures. In Section \ref{pca.sec}, the spiked Wishart model is discussed in detail and generic risk lower and upper bounds are obtained.  The general results are applied in Section \ref{example.sec} to specific settings  and minimax optimal rates are established by  explicitly calculating the local and global metric-entropic quantities. In Section \ref{dis.sec}, we address the computational issues of the proposed estimators  and discuss some extensions and make connections to some other interesting problems. 

For a vector $\bold{a} = (a_1,...,a_n)^\top \in \mathbb{R}^{n}$, we denote $\text{diag}(a_1,...,a_n)\in\R^{n\times n}$ as the diagonal matrix whose $i$-th diagonal entry is $a_i$, and define the $\ell_p$ norm $\| \bold{a} \|_p = \big(\sum_{i=1}^n a_i^p\big)^{1/p}$.  We write $a\land b=\min\{a,b\}$ and $a\lor b=\max\{a,b\}$. For a matrix $ \bold{A}=(a_{ij})\in \R^{p_1\times p_2}$,  we define its Frobenius norm as $\| \bold{A}\|_F = \sqrt{ \sum_{i=1}^{p_1}\sum_{j=1}^{p_2} a^2_{ij}}$ and its spectral norm as $\| \bold{A} \| =\sup_{\|\bold{x}\|_2\le 1}\|\bold{A}\bold{x}\|_2 $; we also denote
$ \bold{A}_{.i}\in \R^{p_1}$ as its $i$-th column and $ \bold{A}_{i.}\in \R^{p_2}$ as its $i$-th row. Let $O(p,k)=\{ \bV\in \R^{p\times k}: \bV^\top \bV={\bf I}_k \}$ be the set of all $p\times k$ orthonormal matrices and $O_p=O(p,p)$, the set of $p$-dimensional orthonormal matrices. For a rank $r$ matrix $\bold{A}\in \R^{p_1\times p_2}$ with $1\le r\le p_1\land p_2$,  its SVD is denoted as $\bold{A}=\bU\bGam\bV^\top$ where $\bU\in O(p_1,r)$, $\bV\in O(p_2,r)$, and $\bGam=\text{diag}(\lambda_1(\bold{A}),\lambda_2(\bold{A}),...,\lambda_r(\bold{A}))$ with $\lambda_{\max}(\bold{A})=\lambda_1(\bold{A})\ge \lambda_2(\bold{A})\ge...\ge\lambda_{p_1\land p_2}(\bold{A})=\lambda_{\min}(\bold{A})\ge 0$ being the ordered singular values of $\bold{A}$. The columns of $\bU$ and the columns of $\bV$ are the left singular vectors and right singular vectors associated to the non-zero singular values of $\bold{A}$, respectively.
For a given set $S$, we denote its cardinality as $|S|$. For sequences $\{a_n\}$ and $\{b_n\}$, we write $a_n = o(b_n)$ or $a_n\ll b_n$ if $\lim_{n} a_n/b_n =0$, and write $a_n = O(b_n)$, $a_n\lesssim b_n$ or $b_n \gtrsim a_n$ if there exists a constant $C$ such that $a_n \le Cb_n$ for all $n$. We write $a_n\asymp b_n$ if $a_n \lesssim b_n$ and $a_n\gtrsim b_n$. Lastly, $c,C, C_0, C_1,...$ are constants that may vary from place to place.

\section{Minimax Lower Bounds via Local Packing} \label{mdm.minimax.sec}

We start with the matrix denoising model. Without loss of generality,  we focus on estimating the structured left singular subspace $\text{span}(\bU)$. Specifically, for a given subset $\CC\subset O(p_1,r)$, we consider the parameter space 
\beq \label{para}
\mathcal{Y}(\CC,t,p_1,p_2,r)=\bigg\{(\bGam,\bU,\bV): \begin{aligned} &\bGam=\text{diag}(\lambda_1,...,\lambda_r),\bU\in \CC, \bV\in O(p_2,r) \\
&  L t\ge \lambda_1\ge...\ge\lambda_r\ge t/L>0
\end{aligned}\bigg\},
\eeq
for some fixed constant $L>1$.
For any $\bU\in O(p_1,r)$ and $\epsilon\in(0,1)$,  the $\epsilon$-ball centered at $\bU$ is defined as
\[
\bB(\bU,\epsilon) = \{ \bU'\in O(p_1,r): d(\bU',\bU)\le\epsilon\},
\]
and for any given subset $\CC\subset O(p_1,r)$, we define
\[
\text{diam}(\CC) = \sup_{\bU_1,\bU_2\in \CC}d(\bU_1,\bU_2).
\]
We introduce the concepts of packing and covering of a given set before stating a general minimax lower bound. 

\begin{definition}[$\epsilon$-packing and $\epsilon$-covering]
	Let $(V,d)$ be a metric space and $M\subset V$. We say that $G(M,d,\epsilon)\subset M$ is an \emph{$\epsilon$-packing} of $M$ if for any $m_i,m_j\in G(M,d,\epsilon)$ with $m_i\ne m_j$, it holds that $d(m_i,m_j)>\epsilon$. We say that $H(M,d,\epsilon)\subset M$ is an  \emph{$\epsilon$-covering} of $M$ if for any $m\in M$, there exists an $m'\in H(M,d,\epsilon)$ such that $d(m,m')<\epsilon$. We denote $\mathcal{M}(M,d,\epsilon)=\max\{|G(M,d,\epsilon)|\}$ and $\mathcal{N}(M,d,\epsilon)=\min \{| H(M,d,\epsilon)|\}$ as the \emph{$\epsilon$-packing number} and the \emph{$\epsilon$-covering number} of $M$, respectively.
\end{definition}

Following \cite{yang1999information}, we also define the metric entropy of a given set.

\begin{definition}[packing and covering $\epsilon$-entropy]
	Let $\mathcal{M}(M,d,\epsilon)$ and $\mathcal{N}(M,d,\epsilon)$ be the $\epsilon$-packing and $\epsilon$-covering number of $M$, respectively. We call $\log \mathcal{M}(M,d,\epsilon)$ the packing $\epsilon$-entropy and  $\log \mathcal{N}(M,d,\epsilon)$ the covering $\epsilon$-entropy of $M$.
\end{definition}

The following theorem gives a minimax lower bound for estimating $\text{span}(\bU)$ over $\mathcal{Y}(\CC,t, p_1,p_2,r)$, as a function of the cardinality of a local packing set of $\CC$, the magnitude of the leading singular values ($t$), the noise level ($\sigma^2$), the rank ($r$), and the dimension ($p_2$) of the right singular vectors in $\bV$. 

\bet \label{mdm.lower.thm}
Under the matrix denoising model $\bY=\bU\bGam\bV^\top+\bZ$ where $(\bGam,\bU,\bV)\in\mathcal{Y}(\CC,t, p_1,p_2,r)$, suppose there exist some $\bU_0\in \CC$, $\epsilon_0>0$ and $\alpha\in(0,1)$ such that a local packing set $G(\bB(\bU_0,\epsilon_0)\cap \CC,d,\alpha\epsilon_0)$ satisfies
\beq\label{cond.lb}
\epsilon_0= \frac{\sqrt{c\sigma^2(t^2+\sigma^2p_2)}}{t^2}\sqrt{\log |G(\bB(\bU_0,\epsilon_0)\cap \CC,d,\alpha\epsilon_0)|}\land \textup{diam}(\CC)
\eeq
for some $c\in(0,1/640]$. Then, as long as $ |G(\bB(\bU_0,\epsilon_0)\cap \CC,d,\alpha\epsilon_0)|\ge 2$, it holds that, for $\theta=(\bGam,\bU,\bV)$,
\beq \label{lb.eq}
\inf_{\widehat{\bU}}\sup_{\theta\in \mathcal{Y}(\CC,t,p_1,p_2,r)}\mathcal{R}(\widehat{\bU},\bU) \gtrsim \bigg( \frac{\sigma\sqrt{t^2+\sigma^2p_2}}{t^2}\sqrt{\log  |G(\bB(\bU_0,\epsilon_0)\cap \CC,d,\alpha\epsilon_0)|}\land \textup{diam}(\CC)\bigg),
\eeq
where the infimum is over all the estimators based on the observation $\bY$.
\eet

The above theorem, to the best of our knowledge, is the first minimax lower bound result for the matrix denoising model under the general parameter space  (\ref{para}).
	Its proof is separated into two parts.  In the strong signal regime ($t^2\gtrsim \sigma^2p_2$), the minimax lower bound can be obtained by generalizing the ideas in \cite{vu2012minimax,vu2013minimax} and \cite{cai2013sparse}, where a general lower bound for testing multiple hypotheses (Lemma \ref{lower.lem})  is applied to obtain (\ref{lb.eq}).  In contrast, the analysis is much more complicated in the weak signal regime ($t^2\lesssim \sigma^2p_2$) due to the asymmetry between $\bU$ and $\bV$: the dependence on $p_2$ need to be captured by extra efforts in the lower bound construction \citep{cai2018rate}, which is different from the aforementioned works on sparse PCA. To achieve this, our analysis relies on a generalized Fano's method for testing multiple composite hypotheses (Lemma \ref{fuzzy.lem.2}) and a nontrivial calculation of the pairwise KL divergence between certain mixture probability measures (Lemma \ref{kl.mix.prop}).

A key observation from the above theorem is the role of the local packing set $G(\bB(\bU_0,\epsilon_0) \cap \CC,d,\alpha\epsilon_0)$ and its entropy measure $\log  |G(\bB(\bU_0,\epsilon_0)\cap \CC,d,\alpha\epsilon_0)|$ in characterizing the fundamental difficulty of the estimation problem. Similar phenomena connecting the local packing numbers to the minimax lower bounds has been observed in, for example, nonparametric function estimation \citep{yang1999information}, high-dimensional linear regression \citep{raskutti2011minimax,verzelen2012minimax}, and sparse principal component analysis \citep{vu2012minimax,cai2013sparse}. 

By \cite{cai2018rate}, a sharp minimax lower bound for estimating $\text{span}(\bU)$ under the unstructured matrix denoising models (i.e., $\CC=O(p_1,r)$) is 
\beq \label{eq.lb.1-5}
\inf_{\widehat{\bU}}\sup_{(\bGam,\bU,\bV)\in \mathcal{Y}(O(p_1,r),t,p_1,p_2,r)}\mathcal{R}(\widehat{\bU},\bU) \gtrsim \bigg( \frac{\sigma\sqrt{(t^2+\sigma^2p_2)rp_1}}{t^2}\land \sqrt{r}\bigg),
\eeq
which, in light of the packing number estimates for the orthogonal group (Lemma 1 of \cite{cai2013sparse}), is a direct consequence of our lower bound (\ref{lb.eq}) for any $\bU_0\in O(p_1,r)$. In addition, comparing the lower bounds (\ref{lb.eq}) and (\ref{eq.lb.1-5}), we observe that the information-geometric quantity ${\log |G(\bB(\bU_0,\epsilon_0)\cap \CC,d,\alpha\epsilon_0)|}$ essentially quantifies the \emph{intrinsic statistical dimension}  (which is $rp_1$ in the case of $\CC=O(p_1,r)$) of the set $\CC$.

\section{Risk Upper Bound using Dudley's Entropy Integral} \label{mle.sec}

In this section, we consider a general singular subspace estimator and study its theoretical properties. Specifically,  we obtain its risk upper bound which, analogous to the minimax lower bound, can be expressed as a function of certain entropic measures related to the structural constraint $\CC$. 

Under the matrix denoising model, with the parameters $(\bGam,\bU,\bV)\in\mathcal{Y}(\CC,t,p_1,p_2,r)$ for some given set $\CC\subset O(p_1,r)$, we consider the structured singular subspace estimator
\beq \label{mle1}
\widehat{\bU}=\argmax_{\bU\in \CC} \text{tr}(\bU^\top\bY\bY^\top \bU).
\eeq
Before stating our main theorem, we need to make more definitions about quantities that play important roles in our subsequent discussions.
\begin{definition}
	For given $\CC\subset O(p_1,r)$ and any $\bU\in \CC$, we define the set
	\[
	\mathcal{T}(\CC,\bU)=\bigg\{\frac{\bW\bW^\top-\bU\bU^\top }{\| \bW\bW^\top-\bU\bU^\top \|_F}\in \R^{p_1\times p_1}: \bW\in\CC\setminus \{\bU\} \bigg\},
	\]
	equipped with the Frobenius distance $d_2$, where for any $\bD_1,\bD_2\in \mathcal{T}(\CC,\bU)$, we define $d_2(\bD_1,\bD_2)=\|\bD_1-\bD_2\|_F$.
\end{definition}
\begin{definition}[Dudley's entropy integral]
	For a metric space $(T,d)$ and a subset $A\subset T$, Dudley's entropy integral of $A$ is defined as $D(A,d)=\int_0^\infty\sqrt{\log \mathcal{N}(A,d,\epsilon)} d\epsilon$. Moreover, we define $D'(A,d)=\int_0^\infty{\log \mathcal{N}(A,d,\epsilon)} d\epsilon$.
\end{definition}
\bet\label{mdm.risk.thm}
Under the matrix denoising model, for any given subset $\CC\subset O(p_1,r)$ and the parameter space $\mathcal{Y}(\CC,t,p_1,p_2,r)$, if $t^2/\sigma^2\gtrsim \sup_{\bU\in\CC} [D'^2(\mathcal{T}(\CC,\bU),d_2)/D^2(\mathcal{T}(\CC,\bU),d_2)]$, it holds that
\beq \label{mdm.risk}
\sup_{(\bGam,\bU,\bV)\in\mathcal{Y}(\CC,t,p_1,p_2,r)}\mathcal{R}(\widehat{\bU},\bU)\lesssim \bigg( \frac{\sigma\Delta(\CC)\sqrt{t^2+\sigma^2p_2}}{t^2} \land \textup{diam}(\CC)\bigg),
\eeq
where $\Delta(\CC)=\sup_{\bU\in \CC}D(\mathcal{T}(\CC,\bU),d_2)$.
\eet

	The proof of the above theorem, as it concerns the generic estimator (\ref{mle1}) under some arbitrary structural set $\CC$, is involved and very different from the existing works such as \cite{cai2013sparse} \cite{deshpande2014cone} \cite{cai2018rate} and \cite{zhang2018heteroskedastic} where specific examples of $\CC$ are considered. The argument relies on careful analysis the supremum of a Gaussian chaos of order 2 and the supremum of a Gaussian process. In the latter case, we applied Dudley's integral inequality (Theorem \ref{dudley.thm}) and the invariance property of the covering numbers with respect to Lipschitz maps (Lemma \ref{szarek.lem}), whereas in the former case, the Arcones-Gin\'e decoupling inequality (Theorem \ref{decouple.thm}) as well as the generic chaining argument (Theorem \ref{chaining.lem}) were used to obtain the desired upper bounds. Many technical tools concatenated for the proof of this theorem can be of independent interest. See more details in Section \ref{proof.upper}.

\begin{figure}[h!]
	\centering
	\includegraphics[angle=0,width=7cm]{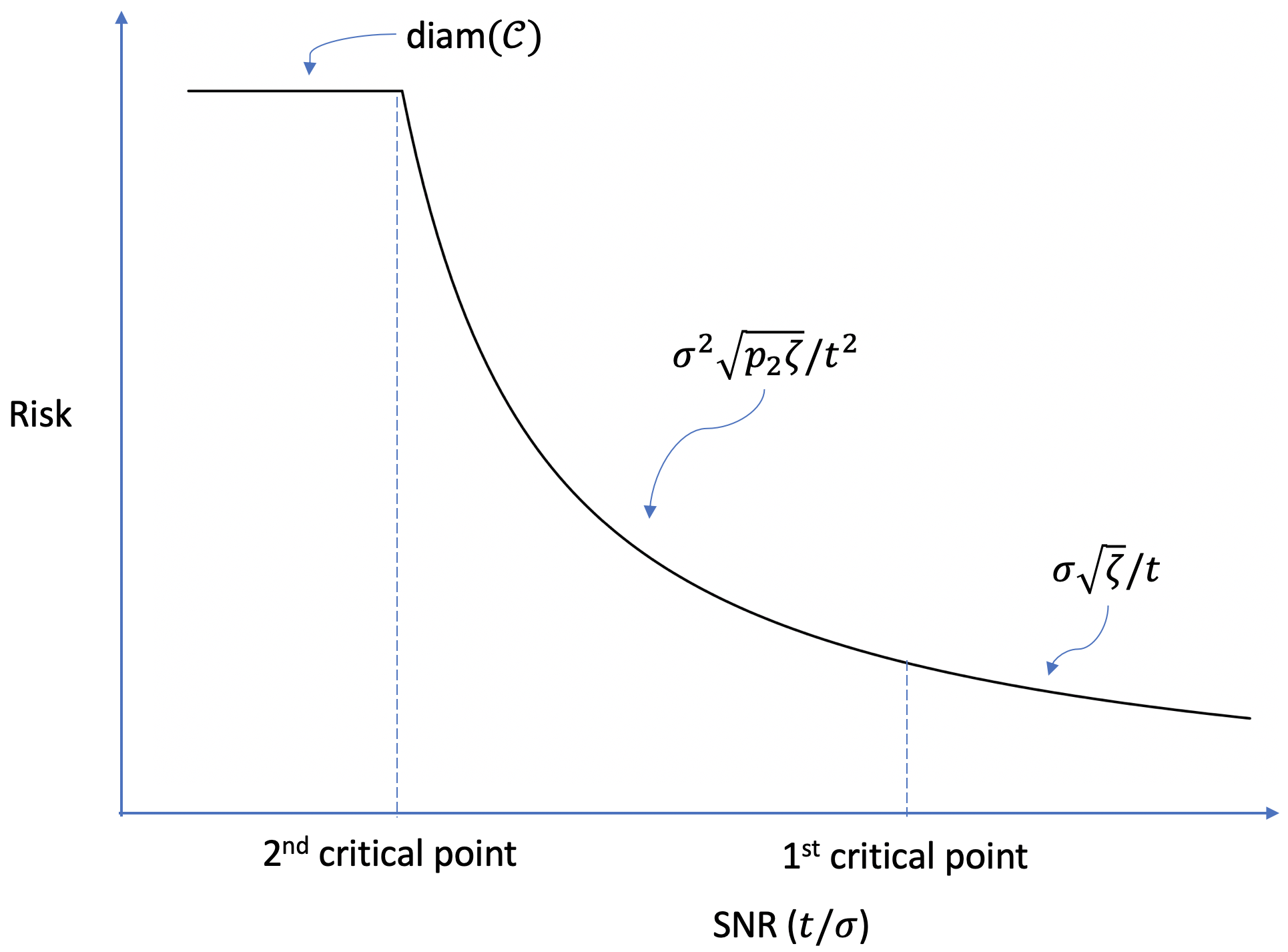}
	\includegraphics[angle=0,width=7cm]{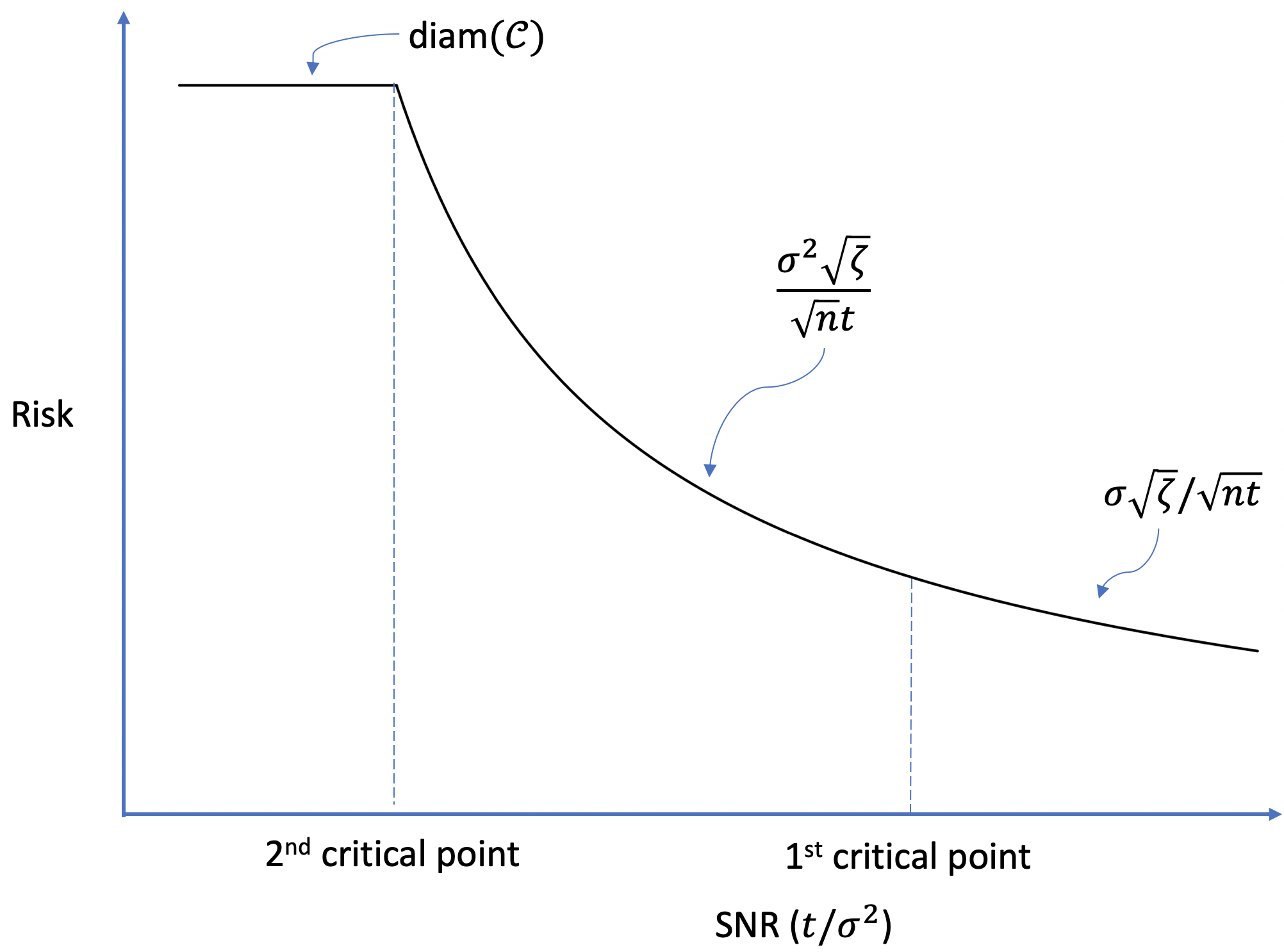}
	\caption{A graphical illustration of the phase transitions in risks as a function of the SNRs under the matrix denoising model (left) and the spiked Wishart model (right).} 
	\label{rate}
\end{figure}

	Interestingly, both the risk upper bound (\ref{mdm.risk}) and the minimax lower bound (\ref{lb.eq}) indicate two phase transitions when treated as a function of the SNR $t/\sigma$, with the first critical point 
	\beq
	\frac{t}{\sigma}\asymp \sqrt{p_2},
	\eeq
	and the second critical point
	\beq \label{cp1}
	\frac{t}{\sigma}\asymp \bigg[\frac{\zeta}{\text{diam}^2(\CC)}+\sqrt{\frac{\zeta p_2}{\text{diam}^2(\CC)}}\bigg]^{1/2},
	\eeq
	where in the upper bound $\zeta=\Delta^2(\CC)$ and in the lower bound  $\zeta=\log  |G(\bB(\bU_0,\epsilon_0)\cap \CC,d,\alpha\epsilon_0)|$. Specifically, the phase transition at the first critical point highlights the role of the dimensionality of the right singular vectors ($\bV$) and the change of the rates of convergence from an inverse quadratic function ($\sigma^2\sqrt{p_2\zeta}/t^2$) to an inverse linear function ($\sigma\sqrt{\zeta}/t$) of $t/\sigma$. The message from the second phase transition concerns the statistical limit of the estimation problem: consistent estimation is possible only when the SNR exceeds the critical point (\ref{cp1}) asymptotically. See Figure \ref{rate} (left) for a graphical illustration. As for the implications of the condition 
	\beq \label{cond}
	t^2/\sigma^2\gtrsim \sup_{\bU\in\CC} [D'^2(\mathcal{T}(\CC,\bU),d_2)/D^2(\mathcal{T}(\CC,\bU),d_2)]
	\eeq
	required by Theorem \ref{mdm.risk.thm}, it can be seen in Section \ref{example.sec} that, for many specific problems, a sufficient condition for (\ref{cond}) is that $t/\sigma$ is above the second critical point (\ref{cp1}), which is mild and natural since the latter condition characterizes the region where $\widehat{\bU}$ is consistent and more generally where consistent estimation is possible.

Comparing our risk upper bound (\ref{mdm.risk}) to the minimax lower bound (\ref{lb.eq}), we can observe the similar role played by the information-geometric quantities that characterize the intrinsic statistical dimension of the sets $\CC$ or $\mathcal{T}(\CC,\bU)$. Specifically, in (\ref{mdm.risk}), the  quantity $\Delta(\CC)$ is related to the global covering entropy, whereas in (\ref{lb.eq}), the quantity $\sqrt{\log  |G(\bB(\bU_0,\epsilon_0)\cap \CC,d,\alpha\epsilon_0)|}$ is associated to the local packing entropy. To obtain the minimax optimal rate of convergence, we need to compare the above two quantities and show 
\beq \label{minimax.rate}
\Delta^2(\CC)\asymp \log  |G(\bB(\bU_0,\epsilon_0)\cap \CC,d,\alpha\epsilon_0)|.
\eeq
Proving the above equation in its general form is difficult. Alternatively, we briefly discuss the affinity between these two geometric quantities yielded by information theory and leave more detailed discussions in the context of some specific examples in Section \ref{example.sec}. 

By definition of the packing numbers, we have the relationship
\beq \label{packing.ineq}
\log  |G(\bB(\bU_0,\epsilon_0)\cap \CC,d,\alpha\epsilon_0)|\le \log \mathcal{M}(\bB(\bU_0,\epsilon_0)\cap \CC,d,\alpha\epsilon_0),
\eeq
that links $\log  |G(\bB(\bU_0,\epsilon_0)\cap \CC,d,\alpha\epsilon_0)|$ to the local packing entropy.
A well known fact about the equivalence between the packing and the covering number of a set $M$ is that
\beq\label{cover.packing}
\mathcal{M}(M,d,2\epsilon)\le \mathcal{N}(M,d,\epsilon)\le \mathcal{M}(M,d,\epsilon).
\eeq
Moreover, \cite{yang1999information} obtained a very interesting result connecting the local and the global (covering) metric entropies. Specifically, let $\bU$ be any element from $M$, then
\beq\label{local.global.entropy}
\log  \mathcal{M}(M,d,\epsilon/2)-\log \mathcal{M}(M,d,\epsilon)\le \log \mathcal{M}(\bB(\bU,\epsilon)\cap M,d,\epsilon/2)\le \log \mathcal{M}(M,d,\epsilon).
\eeq
In Section \ref{example.sec}, by focusing on some specific examples of $\CC$ that are widely considered in practice, we show that  equation (\ref{minimax.rate}) holds, which along with our generic lower and upper bounds recovers some existing minimax rates, and more importantly, helps to establish some previously unknown rates.

\section{Structured Eigen Subspace Estimation in the Spiked Wishart Model} \label{pca.sec}

We turn the focus in this section to the spiked Wishart model where one has i.i.d. observations $Y_i\sim N(\bmu,\bSig)$ with $\bSig=\bU\bGam\bV^\top+\sigma^2 {\bf I}$, which is usually referred as the spiked covariance.
Similar to the matrix denoising model, a minimax lower bound based on some local packing set and a risk upper bound based on the Dudley's entropy integral can be obtained.

\subsection{Minimax Lower Bound}\label{pca.lower.sec}

For any given subset $\CC\subset O(p,r)$, we consider the parameter space 
\[
\mathcal{Z}(\CC,t,p,r)=\{(\bGam,\bU): 
\bGam=\text{diag}(\lambda_1,...,\lambda_r), L t\ge \lambda_1\ge...\ge\lambda_r\ge t/L>0,\bU\in\CC\},
\]
where $L>1$ is some fixed constant.
The following theorem provides minimax lower bound for estimating $\text{span}(\bU)$ over $\mathcal{Z}(\CC,t,p,r)$ under the spiked Wishart model.

\bet \label{lower.bnd.thm.2}
Under the spiked Wishart model where $(\bGam,\bU)\in\mathcal{Z}(\CC,t,p,r)$, suppose there exist some $\bU_0\in\CC$, $\epsilon_0>0$ and $\alpha\in(0,1)$  such that a local packing set $G(\bB(\bU_0,\epsilon_0)\cap \CC,d,\alpha\epsilon_0)$ satisfies
\beq \label{cond.lb2}
\epsilon_0=\frac{\sigma\sqrt{c(\sigma^2+t)}}{t\sqrt{n}}\sqrt{\log |G(\bB(\bU_0,\epsilon_0)\cap \CC,d,\alpha\epsilon_0)|} \land \textup{diam}( \CC),
\eeq
for some $c\in(0,1/32]$. Then, as long as $|G(\bB(\bU_0,\epsilon_0)\cap \CC,d,\alpha\epsilon_0)| \ge 2$, it holds that
\beq \label{lb.eq.2}
\inf_{\widehat{\bU}}\sup_{(\bGam,\bU)\in \mathcal{Z}(\CC,t,p,r)}\mathcal{R}(\widehat{\bU},\bU)  \gtrsim \bigg(\frac{{\sigma} \sqrt{\sigma^2+t}}{t\sqrt{n}}\sqrt{\log |G(\bB(\bU_0,\epsilon_0)\cap \CC,d,\alpha\epsilon_0)| }\land  \textup{diam}( \CC)\bigg),
\eeq
where the infimum is over all the estimators based on the observation $\bY$.
\eet

In \cite{zhang2018heteroskedastic},  a sharp minimax lower bound for estimating $\text{span}(\bU)$ under the unstructured spiked Wishart model was obtained as
\beq \label{lb.eq.3}
\inf_{\widehat{\bU}}\sup_{(\bGam,\bU)\in \mathcal{Z}(O(p,r),t,p,r)}\mathcal{R}(\widehat{\bU},\bU)  \gtrsim \bigg(\frac{{\sigma} \sqrt{(\sigma^2+t)rp}}{t\sqrt{n}}\land  \sqrt{r}\bigg).
\eeq
Comparing the general lower bound (\ref{lb.eq.2}) with (\ref{lb.eq.3}), we observe that the local entropic quantity ${\log  |G(\bB(\bU_0,\epsilon_0)\cap \CC,d,\alpha\epsilon_0)|}$ again characterizes the intrinsic statistical dimension (which is ${rp}$ in the case of $\CC=O(p,r)$) of the set $\CC$. See Section \ref{example.sec} for more examples.

\subsection{Risk Upper Bound} \label{pca.mle.sec}

Under the spiked Wishart model, to estimate the eigen subspace $\text{span}(\bU)$ under the structural constraint $\bU\in \CC$, we start with the sample covariance matrix
\[
\hat{\bSig}=\frac{1}{n}\sum_{i=1}^n(Y_i-\bar{Y})(Y_i-\bar{Y})^\top,
\]
where $\bar{Y}=\frac{1}{n}\sum_{i=1}^nY_i$ and $Y_i$ is the $i$-th row of the observed data matrix $\bY\in \R^{n\times p}$.
Since $\hat{\bSig}$ is invariant to any translation on $\bY$, we assume $\bmu=0$ without loss of generality.

Similar to the matrix denoising model, for the spiked Wishart model, with a slight abuse of notation,  we define the eigen subspace estimator as
\beq \label{mle2}
\widehat{\bU}=\argmax_{\bU\in \CC} \text{tr}(\bU^\top\hat{\bSig}\bU).
\eeq
The following theorem provides the risk upper bound of  $\widehat{\bU}$.

\bet \label{pca.upper.bnd.thm}
Under the spiked Wishart model, for any given $\CC\subset O(p,r)$ and the parameter space $\mathcal{Z}(\CC,t,p,r)$, suppose $n\gtrsim \max\{\log \frac{t}{\sigma^2}, r\}$ and $t/\sigma^2\gtrsim \sup_{\bU\in\CC}[D'^2(\mathcal{T}(\CC,\bU),d_2)/D^2(\mathcal{T}(\CC,\bU),d_2)]$, then
\[
\sup_{(\bGam,\bU)\in \mathcal{Z}(\CC,t,p,r)}\mathcal{R}(\widehat{\bU},\bU) \lesssim \bigg(\frac{\sigma\Delta(\CC)\sqrt{t+\sigma^2}}{t\sqrt{n}} \land\textup{diam}(\CC)\bigg),
\]
where $\Delta(\CC)$ is defined in Theorem \ref{mdm.risk.thm}.
\eet

Similar to the matrix denoising model, the above risk upper bound has a great affinity to the minimax lower bound (\ref{lb.eq.2}), up to a difference in the information-geometric (metric-entropic) measures of $\CC$, and the sharpness of our results relies on the relative magnitude between the  pair of quantities $\Delta^2(\CC)$ and $\log  |G(\bB(\bU_0,\epsilon_0)\cap \CC,d,\alpha\epsilon_0)|$. 
In addition, phase transitions in the rates of the lower and upper bounds as functions of the SNR $t/\sigma^2$ can be observed with the first critical point at
\beq
\frac{t}{\sigma^2}\asymp 1,
\eeq
and the second critical point at
\beq \label{cp2}
\frac{t}{\sigma^2}\asymp \frac{\zeta}{n\cdot\text{diam}^2(\CC)}+\sqrt{\frac{\zeta}{n\cdot\text{diam}^2(\CC)}},
\eeq
where in the lower bound $\zeta= \log|G(\bB(\bU_0,\epsilon_0)\cap \CC,d,\alpha\epsilon_0)| $ and in the upper bound $\zeta=\Delta^2(\CC)$.  Again, the phase transition at the first critical point reflects the change of the speed of the rates of convergence, whereas the phase transition at the second critical point characterizes the statistical limit of the estimation problem. See Figure \ref{rate} (right) for a graphical illustration. Finally, it will be seen in Section \ref{example.sec} that for many specific problems, the condition $t/\sigma^2\gtrsim \sup_{\bU\in\CC}[D'^2(\mathcal{T}(\CC,\bU),d_2)/D^2(\mathcal{T}(\CC,\bU),d_2)]$ required by Theorem \ref{pca.upper.bnd.thm} is mild and in fact necessary for consistent estimation.

\section{Applications} \label{example.sec}

In the following, building upon the minimax lower bounds and the risk upper bounds established in the previous sections, we obtain minimax rates and fundamental limits for various structural principal subspace estimation problems of broad interest. Specifically, in light of our generic results, we show the asymptotic equivalence of various local and global entropic measures associated to some specific examples of $\CC$. Previous discussions under the general settings such as the phase transition phenomena also apply to each of the examples.

\subsection{Sparse PCA/SVD} \label{sparse.sec}

We start with the sparse PCA/SVD where the columns of $\bU$ are sparse vectors.
Suppose $\CC_S(p,r,k)$ is the $k$-sparse subset of $O(p,r)$ for some $k\le p$, i.e., $\CC_S(p,r,k)=\{ \bU\in O(p,r): \max_{1\le i\le r}\|\bU_{.i}\|_0\le k\}.$
The following proposition concerns some estimates about the local and global entropic quantities associated with the set $\CC_S(p,r,k)$. For simplicity, we denote $\CC_S(k)=\CC_S(p,r,k)$ when there is no confusion.

\bep \label{entropy.sparse}
Under the matrix denoising model where $(\bGam,\bU,\bV)\in\mathcal{Y}(\CC_S(k), t,p_1,p_2,r)$ with $k=o(p_1)$ and $r=O(1)$, there exist some $(\bU_0,\epsilon_0,\alpha)$ and a local packing set $G(\bB(\bU_0,\epsilon_0)\cap \CC_S(p_1,r,k),d,\alpha\epsilon_0)$ satisfying (\ref{cond.lb}) such that 
\[
\log |G(\bB(\bU_0,\epsilon_0)\cap \CC_S(p_1,r,k),d,\alpha\epsilon_0)|\asymp\Delta^2(\CC_S(p_1,k,r))\asymp k\log(ep_1/k)+ {k}.
\]
Similarly, under the spiked Wishart model where $(\bGam,\bV)\in\mathcal{Z}(\CC_S(k),t,p,r)$ with $k=o(p)$ and $r=O(1)$, there exist some $(\bU_0,\epsilon_0,\alpha)$ and a local packing set $G(\bB(\bU_0,\epsilon_0)\cap \CC_S(p,r,k),d,\alpha\epsilon_0)$ satisfying (\ref{cond.lb2}) such that 
\[
\log |G(\bB(\bu_0,\epsilon_0)\cap \CC_S(p,k,r),d,\alpha\epsilon_0)|\asymp \Delta^2(\CC_S(p,k,r))\asymp k\log (ep/k)+{k}.
\] 
\eep

In light of our lower and upper bounds under both the matrix denoising model (Theorem \ref{mdm.lower.thm} and \ref{mdm.risk.thm}) and the spiked Wishart model (Theorem \ref{lower.bnd.thm.2} and \ref{pca.upper.bnd.thm}), with Proposition \ref{entropy.sparse}, we are able to establish sharp minimax rates of convergence for sparse PCA/SVD.

\bet \label{sparse.thm}
Under the matrix denoising model with $\bU\in \CC_S(p_1,r,k)$ where $k=o(p_1)$ and $r=O(1)$, it holds that
\beq \label{mdm.sparse}
\inf_{\widehat{\bU}}\sup_{\mathcal{Y}(\CC_S(k), t,p_1,p_2,r)}\mathcal{R}(\widehat{\bU},\bU) \asymp \bigg( \frac{\sigma\sqrt{t^2+\sigma^2 p_2}}{t^2} \bigg( \sqrt{k\log\frac{ep_1}{k}}+\sqrt{k} \bigg)\land 1 \bigg)
\eeq
where the estimator (\ref{mle1}) is rate-optimal whenever consistent estimation is possible.
Similarly, under the spiked Wishart model with $\bU\in\CC_S(p,r,k)$ where $k=o(p)$ and $r=O(1)$,  if $n\gtrsim \max\{\log \frac{t}{\sigma^2}, r\}$, then
\beq \label{swm.sparse}
\inf_{\widehat{\bU}}\sup_{ \mathcal{Z}(\CC_S(k),t,p,r)}\mathcal{R}(\widehat{\bU},\bU) \asymp \bigg( \frac{\sigma\sqrt{t+\sigma^2}}{t\sqrt{n}} \bigg( \sqrt{k\log\frac{ep}{k}}+\sqrt{k} \bigg)\land 1\bigg),
\eeq
where the estimator (\ref{mle2}) is rate-optimal whenever consistent estimation is possible.
\eet

The minimax rate (\ref{swm.sparse}) for the spiked Wishart model (sparse PCA) recovers the ones obtained by \cite{vu2012minimax} and \cite{cai2013sparse} under either rank-one or finite rank $r$ settings. In contrast, the result (\ref{mdm.sparse}) for the matrix denoising model (sparse SVD), to the best of our knowledge, has not been established.

\subsection{Non-Negative PCA/SVD} \label{nnpca.sec}

We now turn to the non-negative PCA/SVD under either the matrix denoising model (SVD) or the spiked Wishart model (PCA) where $\bU\in \CC_{N}(p,r)=\{ \bU=(u_{ij})\in O(p,r):  u_{ij}\ge 0 \text{ for all $i,j$}\}$. 
The following proposition provides estimates about the local and global entropic quantities related to the set $\CC_{N}(p,r)$.

\bep \label{entropy.NN}
Under the matrix denoising model where $(\bGam,\bU,\bV)\in \mathcal{Y}(\CC_{N}(p_1,r),t,p_1,p_2,r)$ and $r=O(1)$, there exist some $(\bU_0,\epsilon_0,\alpha)$ and a local packing set $G(\bB(\bU_0,\epsilon_0)\cap \CC_{N}(p_1,r),d,\alpha\epsilon_0)$ satisfying (\ref{cond.lb}) such that
\[
\Delta^2(\CC_{N}(p_1,r))\asymp \log |G(\bB(\bU_0,\epsilon_0)\cap \CC_{N}(p_1,r),d,\alpha\epsilon_0)|\asymp p_1.
\]
Similarly, under the spiked Wishart model where $(\bGam,\bU)\in\mathcal{Z}(\CC_{N}(p,r),t,p,r)$ and $r=O(1)$, there exist some $(\bU_0,\epsilon_0,\alpha)$ and a local packing set $G(\bB(\bU_0,\epsilon_0)\cap \CC_{N}(p,r),d,\alpha\epsilon_0)$ satisfying (\ref{cond.lb2}) such that 
\[
\Delta^2(\CC_{N}(p,r))\asymp \log |G(\bB(\bU_0,\epsilon_0)\cap \CC_{N}(p,r),d,\alpha\epsilon_0)|\asymp p.
\]
\eep

Proposition \ref{entropy.NN} enables us to establish sharp minimax rates of convergence for non-negative PCA/SVD using the general lower and upper bounds from the  previous sections.

\bet \label{NN.thm}
Under the matrix denoising model with $\bU\in \CC_{N}(p_1,r)$ where $r=O(1)$, it holds that
\beq \label{mdm.nn}
\inf_{\widehat{\bU}}\sup_{ \mathcal{Y}(\CC_{N}(p_1,r),t,p_1,p_2,r)}\mathcal{R}(\widehat{\bU},\bU) \asymp \frac{\sigma\sqrt{(t^2+\sigma^2 p_2)p_1}}{t^2} \land 1,
\eeq
and the estimator  (\ref{mle1}) is rate-optimal whenever consistent estimation is possible.
Similarly, for the spiked Wishart model with $\bU\in\CC_{N}(p,r)$ where $r=O(1)$, if $n\gtrsim \max\{\log \frac{t}{\sigma^2}, r\}$, then
\beq \label{swm.nn}
\inf_{\widehat{\bU}}\sup_{ \mathcal{Z}(\CC_{N}(p,r),t,p,r)}\mathcal{R}(\widehat{\bU},\bU) \asymp  \frac{\sigma\sqrt{(t+\sigma^2)p}}{t\sqrt{n}} \land 1,
\eeq
where the estimator (\ref{mle2}) is rate-optimal whenever consistent estimation is possible.
\eet

The minimax rates for non-negative PCA/SVD, which were previously unknown, turn out to be the same as the  rates for the ordinary unstructured SVD \citep{cai2018rate}  and PCA \citep{zhang2018heteroskedastic}. This is due to the fact claimed in Proposition \ref{entropy.NN} that, under the finite rank scenarios, as a much smaller subset of $O(p,r)$,  $\CC_{N}(p,r)$  has asymptotically the same geometric complexity as $O(p,r)$.

\begin{remark}
	\cite{deshpande2014cone} considered the rank-one Gaussian Wigner model $\bY=\lambda\bu\bu^\top+\bZ\in \R^{p_1\times p_1}$, which can be treated as a special case of the matrix denoising model. Specifically, it was shown that, for $\widehat{\bu}=\argmax_{\bu\in\CC_{N}(p,1) }\bu^\top\bY\bu$, it holds that
	\[
	\sup_{(\lambda,\bu)\in \mathcal{Z}(\CC_{N}(p,1),t,p,1)}\E [1-|\widehat{\bu}^\top \bu|]\lesssim \frac{\sigma\sqrt{p}}{t}\land 1,
	\]
	which, by the fact that $1-|\widehat{\bu}^\top \bu|\le d(\widehat{\bu},\bu)$, can be implied by our result (see also Section \ref{dis.sec2}). Similar problems were studied in  \cite{montanari2015non} under the setting where $p_1/p_2\to \alpha\in (0,\infty)$. However, their focus is to unveil the asymptotic behavior of $\widehat{\bu}^\top \bu$ as well as the analysis of an approximate message passing algorithm, which is different from ours.
\end{remark}

\subsection{Subspace Constrained PCA/SVD} \label{sub.sec}

In some applications such as network clustering \citep{wang2010flexible,kawale2013constrained,kleindessner2019guarantees}, it is of interest to estimate principal subspaces with certain linear subspace constraints. For example, under the matrix denoising model, for some fixed $A\in \R^{p_1\times (p_1-k)}$ of rank $(p_1-k)$ where $r< k< p_1$,  a $k$-dimensional subspace constraint on the singular subspace $\text{span}(\bU)$ could be $\bU\in \CC_A(p_1,r,k)=\{ \bU\in O(p_1,r): A\bU_{.i}=0, \forall 1\le i\le r \}$. Again, subspace constrained PCA/SVD can be solved based on the general results obtained in the previous sections.

\bep \label{entropy.sub}
For  given $A\in \R^{p_1\times (p_1-k)}$ of rank $(p_1-k)$, under the matrix denoising model where $(\bGam,\bU,\bV)\in \mathcal{Y}(\CC_A(p_1,r,k),t,p_1,p_2,r)$ and $r=O(1)$, there exist some $(\bU_0,\epsilon_0,\alpha)$ and a local packing set $G(\bB(\bU_0,\epsilon_0)\cap \CC_A(p_1,r,k),d,\alpha\epsilon_0)$ satisfying (\ref{cond.lb}) such that
\[
\Delta^2(\CC_A(p_1,r,k))\asymp \log |G(\bB(\bu_0,\epsilon_0)\cap \CC_A(p_1,r,k),d,\alpha\epsilon_0)|\asymp k.
\]
Similarly, for given $B\in \R^{p\times (p-k)}$ of rank $(p-k)$, under the spiked Wishart model with $(\bGam,\bU)\in\mathcal{Z}(\CC_B(p,r,k),t,p,r)$ and $r=O(1)$, there exist some $(\bU_0,\epsilon_0,\alpha)$ and a local packing set $G(\bB(\bU_0,\epsilon_0)\cap \CC_B(p,r,k),d,\alpha\epsilon_0)$  satisfying (\ref{cond.lb2}) such that 
\[
\Delta^2(\CC_B(p,r,k)) \asymp \log |G(\bB(\bU_0,\epsilon_0)\cap \CC_B(p,r,k),d,\alpha\epsilon_0)|\asymp k.
\]
\eep

\bet \label{sub.thm}
Under the matrix denoising model with $\bU\in \CC_A(p_1,r,k)$ where $ r< k<p_1$, $r=O(1)$ and $A\in \R^{p_1\times (p_1-k)}$ is of rank $(p_1-k)$, it holds that
\beq \label{mdm.sub}
\inf_{\widehat{\bU}}\sup_{ \mathcal{Y}(\CC_A(p_1,r,k),t,p_1,p_2,r)}\mathcal{R}(\widehat{\bU},\bU) \asymp \bigg( \frac{\sigma\sqrt{(t^2+\sigma^2 p_2)k}}{t^2} \land 1  \bigg)
\eeq
and the estimator (\ref{mle1}) is rate-optimal whenever consistent estimation is possible. Similarly, under the spiked Wishart model with $\bU\in\CC_B(p,r,k)$, where $r< k<p$, $r=O(1)$  and $B\in \R^{p\times (p-k)}$ is of rank $(p-k)$, if $n\gtrsim \max\{\log \frac{t}{\sigma^2}, r\}$, then
\beq \label{swm.sub}
\inf_{\widehat{\bU}}\sup_{\mathcal{Z}(\CC_B(p,r,k),t,p,r)}\mathcal{R}(\widehat{\bU},\bU) \asymp \bigg( \frac{\sigma\sqrt{(t+\sigma^2)k}}{t\sqrt{n}} \land 1\bigg),
\eeq
where the estimator (\ref{mle2}) is rate-optimal whenever consistent estimation is possible.
\eet

\subsection{Spectral Clustering} \label{sc.sec}

As discussed in Section \ref{related.works}, spectral clustering can be treated as estimation of the structural eigenvector under the rank-one matrix denoising model $\bY=\lambda\bu\bv^\top +\bZ\in \R^{n\times p}$ where $\lambda=\|\bh\|_2^2\|\btheta\|_2^2$ is the global signal strength, $\bu=\bh/\|\bh\|_2\in \CC^n_{\pm}=\{\bu\in \R^n:\|\bu\|_2=1, u_i\in\{\pm n^{-1/2}\}\}$ indicates the group labels, and $\bZ$ has i.i.d. entries from $N(0,\sigma^2)$. As a result, important insights about the clustering problem can be obtained by calculating the entropic quantities related to $\CC^n_{\pm}$ and applying the general results from the previous sections.

\bep \label{entropy.pm}
Under the matrix denoising model where $(\lambda,\bu,\bv)\in \mathcal{Y}(\CC^{n}_{\pm},t,n,p,1)$, it holds that $\Delta^2(\CC_{\pm}^n)\lesssim n$. In addition, if $t^2=C \sigma^2(\sqrt{pn}+n)$ for some constant $C>0$, then there exist some $(\bu_0,\epsilon_0,\alpha)$ and a local packing set $G(\bB(\bu_0,\epsilon_0)\cap \CC^{n}_{\pm},d,\alpha\epsilon_0)$ satisfying (\ref{cond.lb}) such that $\log |G(\bB(\bu_0,\epsilon_0)\cap \CC^{n}_{\pm},d,\alpha\epsilon_0)|\asymp n$. 
\eep

\bet \label{pm.thm}
Under the spectral clustering model defined in Section \ref{related.works}, or equivalently, the matrix denoising model $\bY=\lambda\bu\bv^\top +\bZ\in \R^{n\times p}$ where $\bu\in \CC^{n}_{\pm}$, the estimator $\widehat{\bu}=\argmax_{\bu\in \CC^{n}_{\pm}} \bu^\top \bY\bY^\top \bu$ satisfies
\beq \label{mdm.pm}
\sup_{(\lambda,\bu,\bv)\in \mathcal{Y}(\CC^{n}_{\pm},t,n,p,1)}\mathcal{R}(\widehat{\bu},\bu) \lesssim \bigg( \frac{\sigma\sqrt{(t^2+\sigma^2 p)n}}{t^2} \land 1  \bigg).
\eeq
In addition, if $t^2\lesssim \sigma^2(n+\sqrt{np})$, then
\beq
\inf_{\widehat{\bu}}\sup_{(\lambda,\bu,\bv)\in \mathcal{Y}(\CC^{n}_{\pm},t,n,p,1)}\mathcal{R}(\widehat{\bu},\bu) \gtrsim C
\eeq
for some absolute constant $C>0$.
\eet

Intuitively, the fundamental difficulty for clustering relies on the interplay between the global signal strength $\lambda$, which reflects both the sample size ($n$) and the distance between the two clusters ($\|\btheta\|_2$), the noise level ($\sigma^2$), and the dimensionality ($p$). In particular, the lower bound from the above theorem shows that one needs $\lambda^2 \gtrsim \sigma^2(\sqrt{pn}+n)$ in order to have consistent clustering. Moreover, the risk upper bound implies that, whenever $\lambda^2 \gtrsim \sigma^2(\sqrt{pn}+n)$, the estimator $\widehat{\bu}$ is consistent. Theorem \ref{pm.thm} thus establishes  the fundamental statistical limit for the minimal global signal strength for consistent clustering. 
Similar phenomena have also been observed by \cite{azizyan2013minimax} and \cite{cai2018rate}.

Nevertheless, it should be noted that, despite the fundamental limits for consistent recovery yielded by Theorem \ref{pm.thm}, the estimator $\widehat{\bu}$ is in itself sub-optimal and can be further improved through a variant of Lloyd’s iterations. See \cite{lu2016statistical} and \cite{ndaoud2018sharp} for more details.

\section{Discussions} \label{dis.sec}
In this paper, we studied a collection of structural principal subspace estimation problems in a unified framework by  exploring the  deep connections between the difficulty  for statistical estimation and the geometric complexity of the parameter spaces. Minimax optimal rates of convergence for a collection of structured PCA/SVD problems are established. In this section, we discuss the computational issues of the proposed estimators as well as the extensions and connections to other problems.

\subsection{Computationally Efficient Algorithms and the Iterative Projection Method}

In general, the constrained optimization problems that define the estimators in (\ref{mle1}) and (\ref{mle2}) are computationally intractable. However, in practice, many iterative algorithms have been developed to approximate such estimators.

For example, under the matrix denoising model, given the data matrix $\bY$, the set $\CC$, and an initial estimator $\bU_0\in O(p_1,r)$, an iterative algorithm for the constrained optimization problem 
$\argmax_{\bU\in \CC} \text{tr}(\bU^\top \bY\bY^\top \bU)$ 
can be realized through iterations over the following updates for $t\ge 1$:
\begin{enumerate}
	\item Multiplication: $\bG_t=\bY\bY^\top \bU_t$;
	\item QR factorization: $\bU'_{t+1}\bW_{t+1}=\bG_t$ where $\bU'_{t+1}$ is $p_1\times r$ orthonormal and $\bW_{t+1}$ is $r\times r$ upper triangular;
	\item Projection: $\bU_{t+1}= {\PP_{\CC}(\bU'_{t+1})}$.
\end{enumerate}
Here the projection operator $\PP_{\CC}(\cdot)$ is defined as $\PP_{\CC}(\bU)=\argmin_{\bG\in \CC}d(\bU,\bG).$
The above algorithm generalizes the ideas of the projected power method (see, for example, \cite{boumal2016nonconvex,chen2018projected,onaran2017projected}) and the orthogonal iteration method \citep{golub2012matrix,ma2013sparse}.

The computational efficiency of this iterative algorithm relies on the complexity of the projection operator $\PP_{\CC}$ for a given $\CC$. In the rank-one case (r=1), \cite{ferreira2013projections} pointed out that, whenever the set $\CC$ is an intersection of a convex cone and the unit sphere, the projection operator $\PP_{\CC}(\cdot)$ admits an explicit formula and can be computed efficiently. This class of spherical convex sets includes many of the above examples such as non-negative PCA/SVD and subspace constrained PCA/SVD. The case of spectral clustering, under the rank-one setting, is also straightforward as the projection has a simple expression $\PP_{\CC_{\pm}^n}(\bu)=\text{sgn}(\bu)/\sqrt{n}$ (see \cite{ndaoud2018sharp} and \cite{loffler2019optimality} for more in depth discussions). 
As for sparse PCA/SVD, the computational side of the problem is much more complicated and has been extensively studied in literature \citep{shen2008sparse,d2008optimal,witten2009penalized,journee2010generalized,ma2013sparse,vu2013fantope,yuan2013truncated,deshpande2014information}.

In addition to the iterative projection method discussed above, there are several other computationally efficient algorithms such as convex (semidefinite in particular) relaxations \citep{singer2011angular,deshpande2014cone,bandeira2017tightness} and the approximate message passing algorithms \citep{deshpande2014information,deshpande2014cone,montanari2015non,rangan2012iterative}, that have been considered to solve the structured eigenvector problems. However, the focuses of these algorithms are still rank-one matrices, and it remains to be understood how well these algorithms generalize to the general rank-$r$ cases. We leave further investigations along these directions to future work.

\subsection{Extensions and Future Work} \label{dis.sec2}
As mentioned in Section \ref{related.works}, an important special case of matrix denoising model is the {Gaussian Wigner model} \citep{deshpande2014cone,montanari2015non,perry2018optimality}, where the data matrix $\bY= \bU\bGam\bU^\top+\bZ\in \R^{p\times p}$ is symmetric, and the noise matrix $\bZ$ has i.i.d. entries (up to symmetry) drawn from $N(0,\sigma^2)$.  Consider the parameter space $\mathcal{Z}(\CC,t,p,r)$ defined in Section \ref{pca.lower.sec}. It can be shown that, under  similar conditions to those of Theorem \ref{mdm.lower.thm}, 
\beq
\inf_{\widehat{\bU}}\sup_{(\bGam,\bU)\in \mathcal{Z}(\CC,t,p,r)}\mathcal{R}(\widehat{\bU},\bU) \gtrsim \bigg( \frac{\sigma}{t}\sqrt{\log |G(\bB(\bU_0,\epsilon_0)\cap \CC,d,\alpha\epsilon_0 )|}\land \textup{diam}(\CC)\bigg).
\eeq
Moreover, if we define $\widehat{\bU}=\argmax_{\bU\in \CC} \text{tr}(\bU^\top \bY\bU)$, then its risk upper bound can be obtained as
\beq \label{gwm.risk}
\sup_{(\bGam,\bU)\in \mathcal{Z}(\CC,t,p,r)}\mathcal{R}(\widehat{\bU},\bU)\lesssim \bigg( \frac{\sigma{\Delta(\CC)}}{t} \land\textup{diam}(\CC)\bigg).
\eeq
These general bounds combined with the entropic quantities calculated in Section \ref{example.sec} would yield many other interesting optimality results. For instance, recall that the Gaussian $\mathbb{Z}/2$ synchronization problem can be treated as a rank-one Gaussian Wigner model $\bY=\lambda\bu\bu^\top+\bZ$ where $\bu\in\CC_{\pm}^n$. In this case, we have,  for $t\lesssim \sigma\sqrt{n}$
\[
\inf_{\widehat{\bu}}\sup_{(\lambda,\bu)\in \mathcal{Z}(\CC_{\pm}^n,t,n,1)}\mathcal{R}(\widehat{\bu},\bu)\gtrsim C.
\]
and, for $\widehat{\bu}=\argmax_{\bu\in \CC_{\pm}^n} \bu^\top\bY\bu$,
\[
\sup_{(\lambda,\bu)\in \mathcal{Z}(\CC_{\pm}^n,t,n,1)}\mathcal{R}(\widehat{\bu},\bu)\lesssim \bigg( \frac{\sigma\sqrt{n}}{t} \land 1\bigg).
\]
This implies that, about Gaussian $\mathbb{Z}/2$ synchronization, to have consistent estimation/recovery, one needs $\lambda\gtrsim \sigma \sqrt{n}$, and the estimator $\widehat{\bu}$ is consistent whenever $\lambda\gtrsim \sigma\sqrt{n}$.
These results make interesting connections to the existing works \citep{javanmard2016phase,perry2018optimality} concerning the so-called critical threshold or fundamental limit for the SNRs in $\mathbb{Z}/2$ synchronization problems.

In the present paper, under the matrix denoising model, we only focused on the cases where the prior structural knowledge on the targeted singular subspace $\text{span}(\bU)$ is available.  However, in some applications, structural knowledge on the other singular subspace $\text{span}(\bV)$ can also be available. An interesting question is whether and how much the prior knowledge on $\text{span}(\bV)$  will help in the estimation of $\text{span}(\bU)$.  
Some preliminary thinking suggests that novel phenomena might exist in such settings. For example, in an extreme case, if $\bV$ is completely known a priori, then after a simple transform $\bY\bV=\bU\bGam+\bZ\bV$, estimation of $\text{span}(\bU)$ can be reduced to a Gaussian mean estimation problem, whose minimax rate is clearly independent of the dimension of the columns in $\bV$ and therefore quite different from the rates obtained in this paper. The problem again bears important concrete examples in statistics and machine learning. The present work provides a theoretical foundation for studying these problems.

\appendix

\section{Proof of the Main Theorems}

In this section, we prove Theorems \ref{mdm.lower.thm}, \ref{mdm.risk.thm}, \ref{lower.bnd.thm.2} and \ref{pca.upper.bnd.thm}.

\subsection{Risk Upper Bounds} \label{proof.upper}

This section proves Theorems \ref{mdm.risk.thm} and \ref{pca.upper.bnd.thm}. Throughout, for any $\bX,\bY\in \R^{p_1\times p_2}$, we denote $\langle\bX,\bY \rangle=\text{tr}(\bX^\top \bY)$. We recall Lemma 1 in  \cite{cai2018rate}, which concerns the relationships between different distance measures.

\bel\label{dist.lem}
For $\bH_1,\bH_2\in O(p,r)$, $\| \bH_1\bH_1^\top-\bH_2\bH_2^\top\|_F= \sqrt{2(r-\|\bH_1^\top \bH_2\|_F^2)}$, and 
$
\frac{1}{\sqrt{2}}\| \bH_1\bH_1^\top-\bH_2\bH_2^\top\|_F\le  \inf_{\bO\in O(r)}\|\bH_1-\bH_2\bO\|_F\le \| \bH_1\bH_1^\top-\bH_2\bH_2^\top\|_F.
$
\eel

\paragraph{Proof of Theorem \ref{mdm.risk.thm}.} 
We begin by stating a useful lemma, whose proof is delayed to Section \ref{lem.sec}.

\bel \label{F.inn.lem}
Let $\bU\in O(p_1,r)$, and $\bGam=\textup{diag}(\lambda_1,...,\lambda_r)$. Then for any $\bW\in O(p_1,r)$, we have $\frac{\lambda_r^2}{2}\|\bU\bU^\top-\bW\bW^\top\|_F^2\le \langle \bU \bGam^2\bU^\top, \bU\bU^\top-\bW\bW^\top\rangle \le \frac{\lambda_1^2}{2}\|\bU\bU^\top-\bW\bW^\top\|_F^2.$
\eel

By Lemma \ref{F.inn.lem} and the fact that $\text{tr}(\widehat{\bU}^\top\bY\bY^\top \widehat{\bU})\ge \text{tr}({\bU}^\top\bY\bY^\top{\bU})$, or equivalently $\langle \bY\bY^\top, \bU\bU^\top-\widehat{\bU}\widehat{\bU}^\top\rangle\le 0$, we have
\[
\| \widehat{\bU}\widehat{\bU}^\top-\bU\bU^\top\|_F^2 \le \frac{2}{\lambda_r^2} \langle \bU\bGam^2\bU^\top-\bY\bY^\top, \bU\bU^\top-\widehat{\bU}\widehat{\bU}^\top\rangle.
\]
Since $\bY=\bU\bGam\bV^\top+\bZ$, we have $\bY\bY^\top=\bU\bGam^2\bU^\top+\bZ\bV\bGam\bU^\top+\bU\bGam\bV^\top\bZ^\top+\bZ\bZ^\top$. Thus
\begin{align*}
\| \widehat{\bU}\widehat{\bU}^\top-\bU\bU^\top\|_F^2&\le \frac{2}{\lambda_r^2} \big[\langle \bU\bGam\bV^\top\bZ^\top, \widehat{\bU}\widehat{\bU}^\top-\bU\bU^\top\rangle+\langle \bZ\bV\bGam\bU^\top, \widehat{\bU}\widehat{\bU}^\top-\bU\bU^\top\rangle\\
&\quad +\langle \bZ\bZ^\top, \widehat{\bU}\widehat{\bU}^\top-\bU\bU^\top\rangle\big]\\
&\equiv\frac{2}{\lambda_r^2}(H_1+H_2+H_3).
\end{align*}
For $H_1$, if we set
\beq \label{G_W}
\bG_{\bW} = \frac{\bW\bW^\top-\bU\bU^\top}{\|\bW\bW^\top-\bU\bU^\top\|_F}, \qquad\bW\in O(p_1,r)\setminus \{\bU\},
\eeq
we can write
\begin{align*}
H_1&=\langle \bU\bGam\bV^\top\bZ^\top,\widehat{\bU}\widehat{\bU}^\top-\bU\bU^\top\rangle 
=\| \widehat{\bU}\widehat{\bU}^\top-\bU\bU^\top\|_F\cdot \langle \bU\bGam\bV^\top\bZ^\top, \bG_{\widehat{\bU}} \rangle \\
&\le \| \widehat{\bU}\widehat{\bU}^\top-\bU\bU^\top\|_F \cdot \sup_{\bW\in \CC} \text{tr}(\bZ\bV\bGam\bU^\top \bG_{\bW}).
\end{align*}
Similarly, we have $H_2 \le  \| \widehat{\bU}\widehat{\bU}^\top-\bU\bU^\top\|_F \cdot \sup_{\bW\in \CC} \text{tr}(\bU\bGam\bV^\top\bZ^\top \bG_{\bW}),$ and $H_3 \le \| \widehat{\bU}\widehat{\bU}^\top-\bU\bU^\top\|_F \cdot  \sup_{\bW\in \CC} \text{tr}(\bZ^\top \bG_{\bW}\bZ).$
It then follows that 
\beq \label{eq1}
\| \widehat{\bU}\widehat{\bU}^\top-\bU\bU^\top\|_F\le \frac{2}{\lambda_r^2}\bigg(\sup_{\bW\in \CC} \text{tr}(\bZ\bV\bGam\bU^\top \bG_{\bW})+\sup_{\bW\in \CC} \text{tr}(\bU\bGam\bV^\top\bZ^\top \bG_{\bW})+\sup_{\bW\in \CC}\text{tr}(\bZ^\top \bG_{\bW}\bZ)\bigg).
\eeq
The rest of the proof is separated into three parts. In the first two parts, we obtain upper bounds for the  right-hand side of equation (\ref{eq1}). In the third part, we derive the desired risk upper bound.

\underline{Part I.} For the term $\sup_{\bW\in \CC} \text{tr}(\bZ\bV\bGam\bU^\top \bG_{\bW})$, we have
\begin{align*}
&\sup_{\bW\in \CC} \text{tr}(\bZ\bV\bGam\bU^\top \bG_{\bW})=\sup_{\bW\in \CC} \text{tr}(\bU^\top \bG_{\bW}\bZ\bV\bGam)=\sup_{\bW\in \CC} \sum_{i=1}^r \lambda_i (\bU^\top \bG_{\bW}\bZ\bV)_{ii}\\
&\le \lambda_1 \sup_{\bW\in \CC} \text{tr}(\bV\bU^\top \bG_{\bW}\bZ)\le \lambda_1\sup_{\bG\in \mathcal{T}'(\CC,\bU,\bV)} \langle \bG, \bZ\rangle,
\end{align*}
where we defined $\mathcal{T}'(\CC,\bU,\bV)=\{  \bG_{\bW} \bU\bV^\top\in \R^{p_1\times p_2}: \bW\in \CC\setminus \{\bU\}\}.$ To control the expected suprema of the Gaussian process $\sup_{\bG\in \mathcal{T}'(\CC,\bU,\bV)} \langle \bG, \bZ\rangle$, we use the following Dudley's integral inequality (see, for example, \citealt[pp. 188]{vershynin2018high}).

\bet[Dudley's Integral Inequality] \label{dudley.thm}
Let $\{X_t\}_{t\in T}$ be a Gaussian process, that is, a jointly Gaussian family of centered random variables indexed by $T$, where $T$ is equipped with the canonical distance $d(s,t)=\sqrt{\E(X_s-X_t)^2}$. For some universal constant $L$, we have $\E \sup_{t\in T}X_t\le L\int_0^\infty\sqrt{\log \mathcal{N}(T,d,\epsilon)} d\epsilon.$
\eet

For the Gaussian process $\sup_{\bG\in \mathcal{T}'(\CC,\bU,\bV)} \langle \bG, \bZ\rangle$, the canonical distance defined over the set $\mathcal{T}'(\CC,\bU,\bV)$ can be obtained as follows. For any $\bG_1,\bG_2\in \mathcal{T}(\CC,\bU,\bV)$, the canonical distance between $\bG_1$ and $\bG_2$, by definition, is $\sqrt{\E \langle \bG_1-\bG_2,\bZ\rangle^2} = \|\bG_1-\bG_2\|_F\equiv d_2(\bG_1,\bG_2)$.
Theorem \ref{dudley.thm} yields
\beq \label{dudley.T}
\E\sup_{\bG\in \mathcal{T}'(\CC,\bU,\bV)} \langle \bG, \bZ\rangle\le C\sigma\int_0^\infty \sqrt{\log \mathcal{N}(\mathcal{T}'(\CC,\bU,\bV),d_2,\epsilon)}d\epsilon,
\eeq
for some universal constant $C>0$.
Next, for any $\bG_1,\bG_2\in \mathcal{T}'(\CC,\bU,\bV)$, without loss of generality, if we assume $\bG_1=\bG_{\bW_1}\bU\bV^\top$ and $\bG_2=\bG_{\bW_2}\bU\bV^\top$,
where $\bW_1,\bW_2\in \CC\setminus\{\bU\}$, then it holds that
\begin{align} \label{lip.eq}
&d_2(\bG_1,\bG_2)=\| \bG_1-\bG_2\|_F\le \|\bG_{\bW_1}-\bG_{\bW_2}\|_F\|\bU\|\|\bV\|\\
&\le \|\bG_{\bW_1}-\bG_{\bW_2}\|_F=d_2(\bG_{\bW_1},\bG_{\bW_2}),\nonumber 
\end{align}
where we used the fact that $\|\bH\bG\|_F\le \|\bH\|_F\|\bG\|$. The next lemma, obtained by \cite{szarek1998metric}, concerns the invariance property of the covering numbers with respect to Lipschitz maps.
\bel[\cite{szarek1998metric}] \label{szarek.lem}
Let $(M,d)$ and $(M_1,d_1)$ be metric spaces, $K\subset M$, $\Phi:M\to M_1$, and let $L>0$. If $\Phi$ satisfies $d_1(\Phi(x),\Phi(y)) \le Ld(x,y)$ for $x,y,\in M$, then, for every $\epsilon>0$, we have $\mathcal{N}(\Phi(K),d_1,L\epsilon)\le \mathcal{N}(K,d,\epsilon).$
\eel

Define the set $\mathcal{T}(\CC,\bU) = \{ \bG_{\bW}: \bW\in \CC\setminus\{\bU\}\}$. 
Equation (\ref{lip.eq}) and Lemma \ref{szarek.lem} imply
\beq \label{entropy.TC}
\log \mathcal{N}(\mathcal{T}'(\CC,\bU,\bV),d_2,\epsilon)\le \log \mathcal{N}(\mathcal{T}(\CC,\bU),d_2,\epsilon),
\eeq
which means
\beq \label{part1.eq1}
\sup_{\bW\in \CC} \text{tr}(\bZ\bV\bGam\bU^\top \bG_{\bW})\le C\lambda_1\sigma\int_0^\infty \sqrt{\log \mathcal{N}(\mathcal{T}(\CC,\bU),d_2,\epsilon)}d\epsilon.
\eeq
Applying the same argument to $\sup_{\bW\in \CC} \text{tr}(\bU\bGam\bV^\top\bZ^\top \bG_{\bW})$ leads to
\beq \label{part1.eq2}
\sup_{\bW\in \CC} \text{tr}(\bU\bGam\bV^\top\bZ^\top \bG_{\bW})\le C\lambda_1\sigma\int_0^\infty \sqrt{\log \mathcal{N}(\mathcal{T}(\CC,\bU),d_2,\epsilon)}d\epsilon.
\eeq

\underline{Part II.} To bound $\sup_{\bW\in \CC}\text{tr}(\bZ^\top \bG_{\bW}\bZ)$, note that $\text{tr}(\bZ^\top \bG_{\bW}\bZ)= \text{vec}(\bZ)^\top \bD_{\bW} \text{vec}(\bZ)$,
where $\text{vec}(\bZ)=(Z_{11},...,Z_{p_11},Z_{12},...,Z_{p_12},...,Z_{1p_2},...,Z_{p_1p_2})^\top$, and
\beq \label{D.mat}
\bD_{\bW}=\begin{bmatrix}
	\bG_{\bW} & &\\
	&\ddots &\\
	& &  \bG_{\bW} 
\end{bmatrix}\in \R^{p_1p_2\times p_1p_2},
\eeq
It suffices to control the expected supremum of the following Gaussian chaos of order 2,
\beq \label{eq4}
\sup_{\bD\in \mathcal{P}(\CC,\bU)} \text{vec}(\bZ)^\top \bD \text{vec}(\bZ),
\eeq
where $\cP(\CC,\bU)=\{  \bD_{\bW}\in \R^{p_1p_2\times p_1p_2}: \bW\in \CC\setminus\{\bU\} \}$. 
To analyze the above Gaussian chaos, a powerful tool from empirical process theory is the decoupling technique. In particular, we apply the following decoupling inequality obtained by \cite{arcones1993decoupling}  (see also Theorem 2.5 of \cite{krahmer2014suprema}).

\bet[Arcones-Gen\'e Decoupling Inequality] \label{decouple.thm}
Let $\{g_i\}_{1\le i\le n}$ be a sequence of independent standard Gaussian variables and let $\{g'_i\}_{1\le i\le n}$ be an independent copy of  $\{g_i\}_{1\le i\le n}$. Let $\mathcal{B}$ be a collection of $n\times n$ symmetric matrices. Then for all $p\ge 1$, there exists an absolute constant $C$ such that
\[
\E\sup_{B\in \mathcal{B}} \bigg| \sum_{1\le j\ne k\le n} B_{jk}g_jg_k+\sum_{j=1}^nB_{jj}(g_j^2-1)\bigg|^p\le C^p\E \sup_{B\in \mathcal{B}}\bigg|\sum_{1\le j,k\le n} B_{jk}g_jg'_k \bigg|^p.
\]
\eet

From Theorem \ref{decouple.thm} and the fact that for given $\bW\in \CC\setminus\{\bU\}$ we have $\E\text{vec}(\bZ)^\top \bD_{\bW} \text{vec}(\bZ)=0$,  then
\beq \label{decouple.eq}
\E \sup_{\bD\in \cP(\CC,\bU)}[\text{vec}(\bZ)^\top \bD\text{vec}(\bZ)]\le C\E \sup_{\bD\in \cP(\CC,\bU)}[\text{vec}(\bZ)^\top \bD \text{vec}(\bZ')]
\eeq
where $\bZ'$ is an independent copy of $\bZ$. The upper bound of the right hand side of  (\ref{decouple.eq}) can be obtained by using a generic chaining argument developed by \cite{talagrand2014upper}.
To state the result, we make the following definitions that characterize the complexity of a set in a metric space.

\begin{definition}[admissible sequence]
	Given a set $T$ in the metric space $(S,d)$, an admissible sequence is an increasing sequence $\{ \mathcal{A}_n \}$ of partitions of $T$ such that $|\mathcal{A}_0|=1$ and $|\mathcal{A}_n|\le 2^{2^n}$ for $n\ge 1$.
\end{definition}

\begin{definition}[$\gamma_\alpha(T,d)$]
	Given $\alpha>0$ and a set $T$ in the metric space $(S,d)$, we define $\gamma_\alpha(T,d)=\inf\sup_{t\in T}\sum_{n\ge 0} 2^{n/\alpha}\textup{diam}(A_n(t))$,
	where $A_n(t)$ is the unique element of $\mathcal{A}_n$ which contains $t$ and the infimum is taken over all admissible sequences.
\end{definition}

The following theorem from \cite[pp. 246]{talagrand2014upper} provides an important upper bound of the general decoupled Gaussian chaos of order 2.

\bet[\cite{talagrand2014upper}] \label{chaining.lem}
Let $\bg,\bg'\in \R^n$ be independent standard Gaussian vectors, and $\bQ=\{q_{ij}\}_{1\le i,j\le n}\in \R^{n\times n}$.  Given a set $T\subset \R^{n\times n}$ equipped with two distances $d_\infty(\bQ_1,\bQ_2)=\|\bQ_1-\bQ_2\|$ and $d_2(\bQ_1,\bQ_2)=\|\bQ_1-\bQ_2\|_F$,  
\[
\E \sup_{\bQ\in T} \bg^\top \bQ \bg' \le L(\gamma_1(T,d_{\infty})+\gamma_2(T,d_2)),
\]
for some absolute constant $L\ge 0$.
\eet
A direct consequence of Theorem \ref{chaining.lem} is
\beq \label{chaining.eq}
\E \sup_{\bD\in \cP(\CC,\bU)}[\textup{vec}(\bZ)^\top \bD \textup{vec}(\bZ')]\le C\sigma^2(\gamma_1(\cP(\CC,\bU),d_{\infty})+\gamma_2(\cP(\CC,\bU),d_2)).
\eeq
Our next lemma obtains estimates of the functionals $\gamma_1(\cP(\CC,\bU),d_{\infty})$ and $\gamma_2(\cP(\CC,\bU),d_2)$.

\bel \label{entropy.chaos.lem}
Let $\mathcal{T}(\CC,\bU)=\{\bG_{\bW}\in \R^{p_1\times p_1}: \bW\in\CC\setminus\{ \bU\} \}$
be equipped with distances $d_\infty$ and $d_2$ defined in Theorem \ref{chaining.lem}.
It holds that
\begin{align}
\gamma_1(\cP(\CC,\bU),d_\infty)&\le C_1\int_0^{\infty}{\log\mathcal{N}(\mathcal{T}(\CC,\bU),d_2,\epsilon)}d\epsilon,\\
\gamma_2(\cP(\CC,\bU),d_2)&\le C_2\sqrt{p_2}\int_0^{\infty}\sqrt{\log\mathcal{N}(\mathcal{T}(\CC,\bU),d_2,\epsilon)}d\epsilon.
\end{align}
\eel

Combining the above results, we have
\beq \label{part2.eq}
\E\sup_{\bW\in \CC}\text{tr}(\bZ^\top \bG_{\bW}\bZ)\lesssim \sigma^2\sqrt{p_2}\int_0^{\infty}\sqrt{\log\mathcal{N}(\mathcal{T}(\CC,\bU),d_2,\epsilon)}d\epsilon+\sigma^2\int_0^{\infty}{\log\mathcal{N}(\mathcal{T}(\CC,\bU),d_2,\epsilon)}d\epsilon.
\eeq

\underline{Part III.} By (\ref{eq1}) (\ref{part1.eq1}) (\ref{part1.eq2}) and (\ref{part2.eq}), we have,  for any $(\bGam,\bU,\bV)\in \mathcal{Y}(\CC,t,p_1,p_2,r)$, whenever $t\gtrsim \sigma D'(\mathcal{T}(\CC,\bU),d_2)/D(\mathcal{T}(\CC,\bU),d_2)$, 
\begin{align*}
\E \|\widehat{\bU}\widehat{\bU}^\top-\bU\bU^\top\|_F&\lesssim \frac{\sigma \lambda_1D(\mathcal{T}(\CC,\bU),d_2)}{\lambda_r^2}+\frac{\sigma^2\sqrt{p_2}D(\mathcal{T}(\CC,\bU),d_2)+\sigma^2D'(\mathcal{T}(\CC,\bU),d_2)}{\lambda_r^2}\\
&\lesssim \frac{\sigma \Delta(\CC)\sqrt{t^2+\sigma^2p_2}}{t^2}.
\end{align*}
The final result then follows by  noticing the trivial upper bound of $\text{diam}(\CC)$.

\paragraph{Proof of Theorem \ref{pca.upper.bnd.thm}.} We first state a useful lemma (Lemma 3 in  \cite{cai2013sparse}).

\bel\label{F.inn.lem2}
Let $\bSig=\sigma^2{\bf I}_p+\bU\bGam\bU^\top$ where $\bU\in O(p,r)$ and $\bGam=\textup{diag}(\lambda_1,...,\lambda_r)$. Then for any $\bW\in O(p,r)$, we have $\frac{\lambda_r}{2}\|\bU\bU^\top-\bW\bW^\top\|_F^2\le \langle \bSig, \bU\bU^\top-\bW\bW^\top\rangle \le \frac{\lambda_1}{2}\|\bU\bU^\top-\bW\bW^\top\|_F^2.$
\eel

Note that $\bY=\bX \bGam^{1/2} \bU^\top+\bZ\in \R^{n\times p}$ where $\bGam^{1/2}=\text{diag}(\lambda_1^{1/2},...,\lambda_r^{1/2})$,  $\bX\in\R^{n\times r}$ has i.i.d. entries from $\sim N(0, 1)$, and $\bZ$ has i.i.d. entries from $N(0,\sigma^2)$. We can write 
\begin{align*}
\hat{\bSig}=\frac{1}{n}\bY^\top \bY-\bar{Y}\bar{Y}^\top&=\frac{1}{n}( \bU\bGam^{1/2}\bX^\top \bX \bGam^{1/2}\bU^\top +\bZ^\top\bX\bGam^{1/2} \bU^\top +\bU\bGam^{1/2}\bX^\top \bZ+\bZ^\top \bZ) \label{hat.bSig}\\
&\quad -( \bU\bGam^{1/2}\bar{X}\bar{X}^\top\bGam^{1/2}\bU^\top  +\bU\bGam^{1/2}\bar{X}\bar{Z}^\top+\bar{Z}\bar{X}^\top\bGam^{1/2}\bU^\top +\bar{Z}\bar{Z}^\top), \nonumber
\end{align*}
where $\bar{X}=\frac{1}{n}\sum_{i=1}^nX_i\in\R^r$ and $\bar{Z}=\frac{1}{n}\sum_{i=1}^nZ_i\in \R^p$.
Now since $\text{tr}(\widehat{\bU}^\top\hat{\bSig} \widehat{\bU})\ge\text{tr}({\bU}^\top\hat{\bSig} {\bU})$, or equivalently $\langle \hat{\bSig}, \bU\bU^\top-\widehat{\bU}\widehat{\bU}^{\top}\rangle\le 0$, we have
\[
\| \widehat{\bU}\widehat{\bU}^\top-\bU\bU^\top\|_F^2 \le \frac{2}{\lambda_r} \langle \bSig-\hat{\bSig}, \bU\bU^\top-\widehat{\bU}\widehat{\bU}^\top\rangle.
\]
Hence,
\begin{align*}
&\| \widehat{\bU}\widehat{\bU}^\top-\bU\bU^\top\|_F^2\\
&\le \frac{2}{\lambda_r} \big[\langle n^{-1}\bZ^\top\bX\bGam^{1/2} \bU^\top, \widehat{\bU}\widehat{\bU}^\top-\bU\bU^\top\rangle+\langle n^{-1}\bU\bGam^{1/2}\bX^\top \bZ, \widehat{\bU}\widehat{\bU}^\top-\bU\bU^\top\rangle\\
&\quad +\langle  n^{-1} \bU\bGam^{1/2}\bX^\top \bX \bGam^{1/2}\bU^\top-\bU\bGam\bU^\top , \widehat{\bU}\widehat{\bU}^\top-\bU\bU^\top\rangle+\langle  n^{-1}\bZ^\top\bZ-{\bf I}_p, \widehat{\bU}\widehat{\bU}^\top-\bU\bU^\top\rangle\\
&\quad - \langle \bU\bGam^{1/2}\bar{X}\bar{X}^\top\bGam^{1/2}\bU^\top , \widehat{\bU}\widehat{\bU}^\top-\bU\bU^\top\rangle- \langle \bU\bGam^{1/2}\bar{X}\bar{Z}^\top , \widehat{\bU}\widehat{\bU}^\top-\bU\bU^\top\rangle\\
&\quad-\langle \bar{Z}\bar{X}^\top\bGam^{1/2}\bU^\top , \widehat{\bU}\widehat{\bU}^\top-\bU\bU^\top\rangle-\langle \bar{Z}\bar{Z}^\top , \widehat{\bU}\widehat{\bU}^\top-\bU\bU^\top\rangle\big]\\
&\equiv\frac{2}{\lambda_r}(H_1+H_2+H_3+H_4-H_5-H_6-H_7-H_8).
\end{align*}
To control $H_1$,  using the same notations in (\ref{G_W}), we have
\[
H_1\le \frac{1}{n}\| \widehat{\bU}\widehat{\bU}^\top-\bU\bU^\top\|_F\cdot \sup_{\bW\in \CC}\text{tr}(\bU\bGam^{1/2}\bX^\top\bZ\bG_{\bW}).
\]
Similarly, it holds that
\begin{align*}
H_2&\le \frac{1}{n}\| \widehat{\bU}\widehat{\bU}^\top-\bU\bU^\top\|_F\cdot \sup_{\bW\in \CC}\text{tr}(\bZ^\top\bX\bGam^{1/2}\bU^\top\bG_{\bW}),\\
H_3&\le\langle  \bGam^{1/2} (n^{-1}\bX^\top\bX-{\bf I}_r)\bGam^{1/2} , \bU^\top\widehat{\bU}\widehat{\bU}^\top\bU-{\bf I}_r\rangle\\
&\le \|\bGam^{1/2} (n^{-1}\bX^\top\bX-{\bf I}_r)\bGam^{1/2}\|\cdot|\text{tr}(\bU^\top\widehat{\bU}\widehat{\bU}^\top\bU-{\bf I}_r)|\\
&\le \frac{\lambda_1}{2}\|n^{-1} \bX^\top\bX-{\bf I}_r\| \|\bU\bU^\top-\widehat{\bU}\widehat{\bU}^\top\|_F^2,\\
H_4&\le  \| \widehat{\bU}\widehat{\bU}^\top-\bU\bU^\top\|_F\cdot \sup_{\bW\in \CC}\text{tr}((n^{-1}\bZ^\top\bZ-{\bf I}_p)\bG_{\bW}),\\
H_5&\le  \| \bGam^{1/2}\bar{X}\bar{X}^\top\bGam^{1/2} \|\cdot |\text{tr}( \bU^\top\widehat{\bU}\widehat{\bU}^\top\bU-{\bf I}_r)|\le \frac{\lambda_1}{2}\| \bar{X}\bar{X}^\top\|\|\bU\bU^\top-\widehat{\bU}\widehat{\bU}^\top\|_F^2,\\
H_6&\le \| \widehat{\bU}\widehat{\bU}^\top-\bU\bU^\top\|_F\cdot \sup_{\bW\in \CC}\text{tr}(\bU\bGam^{1/2}\bar{X}\bar{Z}^\top\bG_{\bW}),\\
H_7&\le \| \widehat{\bU}\widehat{\bU}^\top-\bU\bU^\top\|_F\cdot \sup_{\bW\in \CC}\text{tr}(\bar{Z}^\top\bar{X}\bGam^{1/2}\bU^\top\bG_{\bW}),\\
H_8&\le  \| \widehat{\bU}\widehat{\bU}^\top-\bU\bU^\top\|_F\cdot \sup_{\bW\in \CC}\text{tr}(\bar{Z}\bar{Z}^\top\bG_{\bW}).
\end{align*}
Combining the above inequalities, we have
\begin{align}
&\| \widehat{\bU}\widehat{\bU}^\top-\bU\bU^\top\|_F \nonumber \\
&\le \frac{2}{\lambda_r(1-\frac{\lambda_1}{\lambda_r}\|n^{-1}\bX^\top\bX-{\bf I}_r\|-\frac{\lambda_1}{\lambda_r}\|\bar{X}\bar{X}^\top\|)}\bigg(n^{-1}\sup_{\bW\in \CC}\text{tr}(\bU\bGam^{1/2}\bX^\top\bZ\bG_{\bW}) \nonumber \\
&\quad+n^{-1}\sup_{\bW\in \CC}\text{tr}(\bZ^\top\bX\bGam^{1/2}\bU^\top\bG_{\bW}) +\sup_{\bW\in \CC}\text{tr}((n^{-1}\bZ^\top\bZ-{\bf I}_p)\bG_{\bW})+ \sup_{\bW\in \CC}\text{tr}(\bU\bGam^{1/2}\bar{X}\bar{Z}^\top\bG_{\bW})\nonumber\\
&\quad+\sup_{\bW\in \CC}\text{tr}(\bar{Z}^\top\bar{X}\bGam^{1/2}\bU^\top\bG_{\bW})+ \sup_{\bW\in \CC}\text{tr}(\bar{Z}\bar{Z}^\top\bG_{\bW})\bigg) \label{pca.eq1} 
\end{align}
The rest of the proof is separated into four parts, with the first three parts controlling the right-hand side of the inequality (\ref{pca.eq1}), and the last part deriving the final risk upper bound.

\underline{Part I.} Note that
\begin{align*}
&\sup_{\bW\in \CC}\text{tr}(\bU\bGam^{1/2}\bX^\top\bZ\bG_{\bW}) =\sup_{\bW\in \CC}\text{tr}(\bX^\top\bZ\bG_{\bW}\bU\bGam^{1/2}) \le \lambda_1^{1/2} \sup_{\bW\in \CC}\text{tr}(\bZ\bG_{\bW}\bU\bX^\top/\|\bX\|)\|\bX\|\\
&\le \lambda_1^{1/2}  \sup_{\bG\in \mathcal{T}_0(\CC,\bU,\bX)}\langle \bZ^\top,\bG\rangle\|\bX\|,
\end{align*}
where $\mathcal{T}_0(\CC,\bU,\bX)=\big\{\frac{\bG_{\bW}\bU\bX^\top}{\|\bX\|}: \bW\in\CC\setminus\{\bU\} \big\}$.
By Theorem \ref{dudley.thm}, we have
\[
\E \bigg[\sup_{\bG\in \mathcal{T}_0(\CC,\bU,\bX)}\langle \bZ^\top,\bG\rangle \bigg| \bX \bigg]\le C\sigma\int_0^\infty \sqrt{\log\mathcal{N}(\mathcal{T}_0(\CC,\bU,\bX),d_2,\epsilon)}d\epsilon.
\]
For any $\bG_1,\bG_2\in \mathcal{T}_0(\CC,\bU,\bX)$, without loss of generality, if we assume $\bG_1=\|\bX\|^{-1}\bG_{\bW_1}\bU\bX^\top$ and $\bG_2=\|\bX\|^{-1}\bG_{\bW_2}\bU\bX^\top$ where $\bW_1,\bW_2\in \CC\setminus\{\bU\}$, then 
\beq \label{lip.eq2}
d_2(\bG_1,\bG_2)\le \|\bG_{\bW_1}-\bG_{\bW_2}\|_F\|\bU\|\le \|\bG_{\bW_1}-\bG_{\bW_2}\|_F=d_2(\bG_{\bW_1},\bG_{\bW_2}).
\eeq
Again, recall the set $\mathcal{T}(\CC,\bU)$ defined in the proof of Theorem \ref{mdm.risk.thm}, by Lemma \ref{szarek.lem}, we have $\log\mathcal{N}(\mathcal{T}_0(\CC,\bU,\bX),d_2,\epsilon) \le \log\mathcal{N}(\mathcal{T}(\CC,\bU),d_2,\epsilon)$, which implies
\beq
\E\sup_{\bW\in \CC}\text{tr}(\bU\bGam^{1/2}\bX^\top\bZ\bG_{\bW}) \le  C\lambda_1^{1/2} \E\|\bX\|\sigma\int_0^\infty \sqrt{\log\mathcal{N}(\mathcal{T}(\CC,\bU),d_2,\epsilon)}d\epsilon.
\eeq
Now by Theorem 5.32 of \cite{vershynin2010introduction}, we have $\E \|\bX\|\le \sqrt{n}+\sqrt{r}$, then
\beq \label{pca.eq2}
\E n^{-1}\sup_{\bW\in \CC}\text{tr}(\bU\bGam^{1/2}\bX^\top\bZ\bG_{\bW}) \le   C\lambda_1^{1/2} \sigma(1/\sqrt{n}+\sqrt{r}/n)\int_0^\infty \sqrt{\log\mathcal{N}(\mathcal{T}(\CC,\bU),d_2,\epsilon)}d\epsilon.
\eeq
Similarly, we can derive
\beq \label{pca.eq3}
\E n^{-1}\sup_{\bW\in \CC}\text{tr}(\bZ^\top \bX\bGam^{1/2} \bU^\top\bG_{\bW}) \le   C\lambda_1^{1/2} \sigma(1/\sqrt{n}+\sqrt{r}/n)\int_0^\infty \sqrt{\log\mathcal{N}(\mathcal{T}(\CC,\bU),d_2,\epsilon)}d\epsilon.
\eeq
One the other hand, since
$
\sup_{\bW\in \CC}\text{tr}(\bU\bGam^{1/2}\bar{X}\bar{Z}^\top\bG_{\bW})=\sup_{\bW\in \CC}\text{tr}(\bar{X}\bar{Z}^\top\bG_{\bW}\bU\bGam^{1/2})\le \lambda_1^{1/2} \sup_{\bW\in \CC}\text{tr}(\bar{Z}^\top\bG_{\bW}\bU\bar{X}/\|\bar{X}\|_2)\|\bar{X}\|_2
\le \lambda_1^{1/2} \|\bar{X}\|_2\sup_{\bg\in \mathcal{T}_1}\langle \bar{Z},\bg \rangle,
$
where $\mathcal{T}_1(\CC,\bU,\bX)=\{ \frac{\bG_{\bW}\bU\bar{X}}{\|\bar{X}\|_2}:\bW\in\CC\setminus\{\bU\} \}$ is equipped with the Euclidean $\ell_2$ distance. By Theorem \ref{dudley.thm}, we have
\[
\E \bigg[ \sup_{\bg\in \mathcal{T}_1(\CC,\bU,\bX)}\langle \bar{Z},\bg \rangle\bigg| \bX \bigg]\le \frac{C\sigma}{\sqrt{n}}\int_0^{\infty}\sqrt{\log \mathcal{N}(\mathcal{T}_1(\CC,\bU,\bX),d_2,\epsilon)}d\epsilon.
\]
Now for any $\bg_1,\bg_2\in \mathcal{T}_1(\CC,\bU,\bX)$, without loss of generality, if we assume $\bg_1=\|\bar{X}\|_2^{-1} \bG_{\bW_1}\bU\bar{X}$ and $\bg_2=\|\bar{X}\|_2^{-1} \bG_{\bW_2}\bU\bar{X},$
where $\bW_1,\bW_2\in \CC\setminus\{\bU\}$, then $\|\bg_1-\bg_2\|_2\le  \|\bar{X}\|_2^{-1} \|\bG_{\bW_1}\bU\bar{X}- \bG_{\bW_2}\bU\bar{X}\|_2\le d_{\infty}(\bG_{\bW_1}, \bG_{\bW_2})\le d_2(\bG_{\bW_1}, \bG_{\bW_2}).$
Lemma \ref{szarek.lem} implies $\log \mathcal{N}(\mathcal{T}_1(\CC,\bU,\bX),d_2,\epsilon)\le \log \mathcal{N}(\mathcal{T}(\CC,\bU),d_2,\epsilon)$,
which along with the fact that $\E \|\bar{X}\|_2\lesssim \sqrt{r/n}$ implies
\beq \label{pca.eq5}
\E \sup_{\bW\in \CC}\text{tr}(\bU\bGam^{1/2}\bar{X}\bar{Z}^\top\bG_{\bW})\le \frac{ C\sigma\sqrt{r}\lambda_1^{1/2}}{{n}} \int_0^{\infty}\sqrt{\log \mathcal{N}(\mathcal{T}(\CC,\bU),d_2,\epsilon)}d\epsilon.
\eeq
Similarly, we have
\beq \label{pca.eq6}
\sup_{\bW\in \CC}\text{tr}(\bar{Z}^\top\bar{X}\bGam^{1/2}\bU^\top\bG_{\bW})\le \frac{ C\sigma\sqrt{r}\lambda_1^{1/2}}{{n}} \int_0^{\infty}\sqrt{\log \mathcal{N}(\mathcal{T}(\CC,\bU),d_2,\epsilon)}d\epsilon.
\eeq

\underline{Part II.} Note that $\text{tr}((n^{-1}\bZ^\top\bZ-\sigma^2{\bf I}_p)\bG_{\bW})=\text{tr}(n^{-1}\bZ^\top\bZ\bG_{\bW})-\sigma^2\text{tr}(\bG_{\bW})=n^{-1}\text{vec}(\bZ)^\top \bD_{\bW} \text{vec}(\bZ),$
where $\bD_{\bW}$ is defined in (\ref{D.mat}).
By the similar chaining argument in Part II of the proof of Theorem \ref{mdm.risk.thm}, we have
\beq \label{pca.eq4}
\E \sup_{\bW\in\CC}\text{tr}((n^{-1}\bZ^\top\bZ-{\bf I}_p)\bG_{\bW})\lesssim \frac{\sigma^2}{\sqrt{n}}\int_0^{\infty} \sqrt{\log \mathcal{N}(\mathcal{T}(\CC,\bU),d_2,\epsilon)}d\epsilon+\frac{\sigma^2}{n}\int_0^{\infty}{\log \mathcal{N}(\mathcal{T}(\CC,\bU),d_2,\epsilon)}d\epsilon
\eeq
Similarly, since $\sup_{\bW\in \CC}\text{tr}(\bar{Z}\bar{Z}^\top\bG_{\bW})=\sup_{\bW\in \CC}\bar{Z}^\top\bG_{\bW}\bar{Z}$, we also have
\beq \label{pca.eq7}
\E \sup_{\bW\in \CC}\text{tr}(\bar{Z}\bar{Z}^\top\bG_{\bW})\lesssim \frac{\sigma^2}{n}\int_0^{\infty}\sqrt{\log \mathcal{N}(\mathcal{T}(\CC,\bU),d_2,\epsilon)}d\epsilon+\frac{\sigma^2}{n}\int_0^{\infty}{\log \mathcal{N}(\mathcal{T}(\CC,\bU),d_2,\epsilon)}d\epsilon.
\eeq

\underline{Part III.} Define the event $E=\{ \|n^{-1}\bX^\top\bX-{\bf I}_r\|\le 1/(4L^2), \|\bar{X}\bar{X}^\top\|\le 1/(4L^2) \}$, where $L$ is the constant in $\mathcal{Z}(\CC,t,p,r)$.
By Proposition D.1 in the Supplementary Material of \cite{ma2013sparse}, 
\[
P(\|n^{-1} \bX^\top \bX-{\bf I}_r\| \le 2(\sqrt{r/n}+t)+(\sqrt{r/n}+t)^2)\ge 1-2e^{-nt^2/2},
\]
which implies $P(\|n^{-1} \bX^\top \bX-{\bf I}_r\| \le 1/(4L^2))\ge 1-2e^{-cn}.$
In addition, since $\|\bar{X}\bar{X}^\top\|\le \|\bar{X}\|_2^2=\frac{1}{n}\sum_{i=1}^rg_i^2$,
where $g_i\sim_{i.i.d.} N(0,1)$, it follows from the  concentration inequality for independent exponential variables  \cite[Proposition 5.16]{vershynin2010introduction} that $P(\|\bar{X}\bar{X}^\top\|\le 1/(4L^2))\ge 1-2e^{-cn}$.
Thus, it follows that
\begin{align*}
P(E^c)\le P(\|n^{-1} \bX^\top \bX-{\bf I}_r\| &\ge 1/(4L^2))+P(\|\bar{X}\bar{X}^\top\|\ge 1/(4L^2))\le 4e^{-cn}.
\end{align*}

\underline{Part IV.} Note that
$\E d(\bU,\widehat{\bU})=\E [d(\bU,\widehat{\bU})|E]+\E [d(\bU,\widehat{\bU})|E^c].$
It follows from (\ref{pca.eq1}) and the inequalities (\ref{pca.eq2})-(\ref{pca.eq7}) from Parts I and II that
\begin{align*}
&\sup_{(\bGam,\bU)\in \mathcal{Z}(\CC,t,p,t)}\E [d(\bU,\widehat{\bU})|E]\\
&\le \frac{C}{t}\bigg[\sqrt{t}\sigma\bigg(\frac{1}{\sqrt{n}}+\frac{\sqrt{r}}{n}\bigg) D(\mathcal{T}(\CC,\bU),d_2)+\frac{\sigma^2D(\mathcal{T}(\CC,\bU),d_2)}{\sqrt{n}}+\frac{\sigma^2D'(\mathcal{T}(\CC,\bU),d_2)}{{n}}\bigg]\\
&\le \frac{C \sigma\Delta(\CC)\sqrt{t(1+r/n)+\sigma^2}}{\sqrt{n}t},
\end{align*} 
where the last inequality holds whenever $t/\sigma^2\gtrsim \sup_{\bU\in\CC}[ D'^2(\mathcal{T}(\CC,\bU),d_2)/D^2(\mathcal{T}(\CC,\bU),d_2)]$.
On the other hand, by Part III, $\E[d(\bU,\widehat{\bU})|E^c]\le \text{diam}(\CC)\cdot P(E^c)\le C\sqrt{r} e^{-cn}$.
Consequently as long as $n\gtrsim \max\{\log \frac{t}{\sigma^2}, r\}$ and $t/\sigma^2\gtrsim \sup_{\bU\in\CC}[ D'^2(\mathcal{T}(\CC,\bU),d_2)/D^2(\mathcal{T}(\CC,\bU),d_2)]$, we have
\[
\sup_{(\bGam,\bU)\in \mathcal{Z}(\CC,t,p,t)}\E d(\bU,\widehat{\bU})\le \frac{C \sigma\Delta(\CC)\sqrt{t+\sigma^2}}{\sqrt{n}t }.
\]
The final result then follows by noticing the trivial upper bound of $\text{diam}(\CC)$.

\subsection{Minimax Lower Bounds}  \label{proof.lower}

\paragraph{Proof of Theorem \ref{mdm.lower.thm}.}
The proof is divided into two parts, the strong signal regime ($t^2\ge \sigma^2p_2/4$) and  the weak signal regime ($t^2<\sigma^2p_2/4$).

\underline{Part I. Strong Signal Regime.}
The following general lower bound for testing multiple hypotheses \cite{tsybakov2009introduction} are needed.

\bel[\cite{tsybakov2009introduction}] \label{lower.lem}
Assume that $M\ge 2$ and suppose that $(\Theta,d)$ contains elements $\theta_0,\theta_1,...,\theta_M$ such that: (i) $d(\theta_j,\theta_k)\ge 2s>0$ for any $0\le  j<k\le M$; (ii) it holds that $\frac{1}{M}\sum_{j=1}^M D(P_j,P_0)\le \alpha \log M$ with $0<\alpha<1/8$ and $P_j=P_{\theta_j}$ for $j=0,1,...,M$,  where $D(P_j,P_0)=\int\log \frac{dP_j}{dP_0}dP_j$ is the KL divergence between $P_j$ and $P_0$. Then
\[
\inf_{\hat{\theta}}\sup_{\theta\in\Theta}P_\theta(d(\hat{\theta},\theta)\ge s) \ge \frac{\sqrt{M}}{1+\sqrt{M}}\bigg( 1-2\alpha-\sqrt{\frac{2\alpha}{\log M}}\bigg).
\]
\eel 

Let $\bV_0\in O(p_2,r)$ be fixed and $\bU_0\in \CC$.
Denote the $\epsilon$-ball $B(\bU_0,\epsilon) = \{ \bU\in O(p_1,r): d(\bU,\bU_0)\le\epsilon\}.$
For some $\delta<\epsilon$, we consider the local $\delta$-packing set $G_\delta=G(B(\bU_0,\epsilon)\cap \CC,d,\delta)$ such that for any pair $\bU,\bU'\in B(\bU_0,\epsilon)\cap \CC$, it holds that $d(\bU,\bU')=\|\bU\bU^\top-\bU'\bU'^\top\|_F\ge \delta.$
We denote the elements of $G_\delta$ as $\bU_i$ for $1\le i\le |G_\delta|$. Lemma \ref{dist.lem} shows that, for any $i$, we can find $\bO_i\in O_r$ such that $\|\bU_0-\bU_i\bO_i\|_F\le d(\bU_0,\bU_i)\le\epsilon.$
Set $\bU'_i=\bU_i\bO_i$ and denote $G'_{\delta}=\{ \bU'_i\}$.
For given $t> 0$, we consider the subset
\[
\mathcal{X}(t,\epsilon,\delta,\bU_0,\bV_0)=\{ (\bGam,\bU,\bV): \bU \in G'_\delta,\bV=\bV_0, \bGam=t{\bf I}_r \}\subset \mathcal{Y}(\CC,t,p_1,p_2,r).
\]
In particular, the above construction admits $|\mathcal{X}(t,\epsilon,\delta,\bU_0,\bV_0)|= |G_\delta|.$

Moreover, for any $(\bGam,\bU_i,\bV_0) \in \mathcal{X}(t, \epsilon,\delta,\bU_0,\bV_0)$, let $P_i$ be the probability measure of $\bY=\bU_i\bGam\bV_0^\top+\bZ$ where $\bZ$ has i.i.d. entries from $N(0,\sigma^2)$. We have, for $1\le i\ne j\le |G_\delta|$,
\begin{align*}
D(P_i,P_j) &=\frac{\|(\bU'_i-\bU'_j)\bGam\bV_0^\top\|_F^2}{2\sigma^2}\le \frac{t^2\|\bU'_i-\bU'_j\|_F^2}{2\sigma^2}\le\frac{2t^2\epsilon^2}{\sigma^2}.
\end{align*}
Now set $\epsilon=\epsilon_0$ and $\delta=\alpha\epsilon$ for some $\alpha\in(0,1)$. By assumption,
\beq \label{epsilon0}
\bigg(\frac{c\sigma^2}{t^2}\log |G_{\alpha\epsilon_0}| \land \text{diam}^2(\CC)\bigg)\le \epsilon_0^2\le \bigg(\frac{\sigma^2}{32t^2}\log |G_{\alpha\epsilon_0}| \land \text{diam}^2(\CC)\bigg)
\eeq
for some $c\in(0,1/32)$, it holds that $D(P_i,P_j) \le \frac{1}{16}\log |G_{\alpha\epsilon_0}|.$
Now by Lemma \ref{lower.lem}, it holds that, for $\theta=(\bGam,\bU,\bV)$,
\[
\inf_{\widehat{\bU}}\sup_{\theta\in \mathcal{X}(t,\epsilon,\delta,\bU_0,\bV_0)} P_\theta(d(\widehat{\bU}, \bU )\ge \alpha\epsilon_0/2 )\ge \frac{\sqrt{ |G_{\alpha\epsilon_0}|}}{1+\sqrt{|G_{\alpha\epsilon_0}|}}\bigg( \frac{7}{8}-\frac{1}{\sqrt{8\log  |G_{\alpha\epsilon_0}|}}\bigg).
\]
By Markov's inequality, we have
\begin{align*}
\inf_{\widehat{\bU}}\sup_{\theta\in \mathcal{X}(t,\epsilon,\delta,\bU_0,\bV_0)}\E_\theta d(\widehat{\bU}, \bU )&\ge \frac{\alpha\epsilon_0\sqrt{ |G_{\alpha\epsilon_0}|}}{2(1+\sqrt{ |G_{\alpha\epsilon_0}|})}\bigg( \frac{7}{8}-\frac{1}{\sqrt{8\log  |G_{\alpha\epsilon_0}|}}\bigg)\ge C \alpha\epsilon_0,
\end{align*}
for some $C>0$ as long as $ |G_{\alpha\epsilon_0}|\ge 2$. Therefore, it holds that
\begin{align*}
&\inf_{\widehat{\bU}}\sup_{\theta\in \mathcal{Y}(\CC,t,p_1,p_2,r)}\E_\theta d(\widehat{\bU}, \bU )\ge \inf_{\widehat{\bU}}\sup_{\theta\in \mathcal{X}(t,\epsilon,\delta,\bU_0,\bV_0)}\E_\theta d(\widehat{\bU}, \bU )\\
&\gtrsim ({\sigma} t^{-1}\sqrt{\log |G_{\alpha\epsilon_0}|}\land  \text{diam}( \CC))\gtrsim \bigg(\frac{\sigma\sqrt{t^2+\sigma^2p_2}}{t^2}\sqrt{\log |G_{\alpha\epsilon_0}|}\land  \text{diam}( \CC)\bigg).
\end{align*}

\underline{Part II. Weak Signal Regime.}
The proof relies on the following generalized Fano's method, obtained by \cite{ma2019optimalb}, about testing multiple composite hypotheses.

\bel[Generalized Fano's Method] \label{fuzzy.lem.2}
Let $\mu_0,\mu_1,...,\mu_M$ be $M+1$ priors on the parameter spaces $\Theta$ of the family $\{P_\theta\}$, and let $P_j$ be the posterior probability measures on $(\mathcal{X},\mathcal{A})$ such that
\[
P_j(S)=\int P_\theta(S)\mu_j(d\theta),\quad \forall S\in\mathcal{A},\quad j=0,1,...,M.
\]
Let $F:\Theta\to (\R^d,d)$. If (i) there exist some sets $B_0,B_1,...,B_M\subset \R^d$ such that $d(B_i,B_j)\ge 2s$ for some $s>0$ for all $0\le i\ne j\le M$ and $\mu_j(\theta\in\Theta: F(\theta)\in B_j)= 1$; and (ii) it holds that $\frac{1}{M}\sum_{j=1}^MD(P_j,P_0)\le \alpha \log M$ 	with $0<\alpha<1/8$. Then
\[
\inf_{\hat{F}}\sup_{\theta\in\Theta} P_{\theta}(d(\hat{F},F(\theta))\ge s)\ge \frac{\sqrt{M}}{1+\sqrt{M}}\bigg(1-2\alpha-\sqrt{\frac{2\alpha}{\log M}} \bigg).
\]
\eel

To use the above lemma, we need to construct a collection of priors over the set $\mathcal{Y}(\CC,t,p_1,p_2,r)$. Specifically, recall the previously constructed $\delta$-packing set $G_\delta=\{\bU_i: 1\le i\le|G_{\delta}|\}$. Inspired by \cite{cai2018rate}, we consider the prior probability measure $\mu_i$ over $\mathcal{Y}(\CC,t,p_1,p_2,r)$, whose definition is given as follows. Let $\bW$ be a random matrix on $\R^{p_2\times r}$, whose probability density is given by
\[
p(\bW)= C\bigg( \frac{p_2}{2\pi}\bigg)^{rp_2/2}\exp(-p_2\|\bW\|_F^2/2)\cdot 1\{1/2\le \lambda_{\min}(\bW)\le \lambda_{\max}(\bW)\le  2\},
\]
where $C$ is a normalizing constant; then, if we denote $\tilde{\bU}_i\tilde{\bGam}_i\tilde{\bV}_i^\top$ as the SVD of $t\bU_i\bW^\top\in \R^{p_1\times p_2}$ where $\bU_i\in G_\delta$ and $\bW\sim p(\bW)$, then  $\mu_i$ is defined as the joint distribution of $(\tilde{\bGam}_i,\tilde{\bU}_i,\tilde{\bV}_i)$. By definition of $\bU_i$, one can easily verify that $\mu_i$ is a well-defined probability measure on $\mathcal{Y}(\CC,t,p_1,p_2,r)$. Note that, for any $\theta_i=(\tilde{\bGam}_i,\tilde{\bU}_i,\tilde{\bV}_i)\in \supp(\mu_i)$ and  $\theta_j=(\tilde{\bGam}_j,\tilde{\bU}_j,\tilde{\bV}_j)\in \supp(\mu_j)$ with $1\le i\ne j\le |G_{\delta}|$, it holds that $d(\tilde{\bU}_i,\tilde{\bU}_j)=d({\bU}_i,{\bU}_j)\ge \delta$.

Consequently, the joint distribution of $\bY=\bU\bGam\bV^\top+\bZ$ with $(\bGam,\bU,\bV)\sim \mu_i$ and $Z_{ij}\sim N(0,\sigma^2)$ can be expressed as
\begin{align*}
P_{i}(\bY)=C\int_{\substack{1/2\le \lambda_{\min}(\bW)\le\lambda_{\max}(\bW)\le  2 }} &\frac{\sigma^{-p_1p_2}}{(2\pi)^{p_1p_2/2}}\exp(-\|\bY-t\bU_i\bW^\top\|_F^2/(2\sigma^2))\\
&\times\bigg( \frac{p_2}{2\pi}\bigg)^{rp_2/2}\exp(-p_2\|\bW\|_F^2/2)d\bW,
\end{align*}
and it remains to control the pairwise KL divergence $D(P_i,P_j)$ for any $1\le i\ne j \le |G_{\delta}|.$ This is done by the next lemma, whose proof, which is involved, is delayed to  Section \ref{lem.sec}.

\bel \label{kl.mix.prop}
Under the assumption of the theorem, for any $1\le i\ne j \le |G_{\delta}|$, we have $D({P}_{i}, {P}_{j}) \le \frac{C_1t^4d^2(\bU_i,\bU_j)}{\sigma^2(4t^2+\sigma^2p_2)}+C_{2}$
where $C_{1},C_2>0$ are some uniform constant and $\{\bU_i\}$ are elements of $G_\delta$.
\eel

Again, set $\epsilon=\epsilon_0$ and $\delta=\alpha\epsilon$ for some $\alpha\in (0,1)$. By assumption,
\[
\bigg(\frac{c\sigma^2(t^2+\sigma^2p_2)}{t^4}\log| G_{\delta}|\land \text{diam}(\CC)\bigg)\le \epsilon_0^2\le \bigg(\frac{\sigma^2(t^2+\sigma^2p_2)}{640t^4}\log| G_{\delta}|\land \text{diam}(\CC)\bigg),
\]
for some $c\in(0,1/640]$.
It then follows that $D(P_i,P_j) \le C\log|G_{\alpha\epsilon_0}|+C_2.$
Now let $\mathcal{X}'(t,\epsilon,\delta,\bU_0)=\bigcup_{1\le i\le |G_{\alpha\epsilon_0}|} \supp(\mu_i).$
By Lemma \ref{fuzzy.lem.2} and Markov's inequality, we have, for $\theta=(\bGam,\bU,\bV)$,
\begin{align*}
\inf_{\widehat{\bU}}\sup_{\theta\in \mathcal{X}'(t,\epsilon_0,\alpha\epsilon_0,\bU_0)}\E_\theta d(\widehat{\bU}, \bU )&\ge \frac{\alpha\epsilon_0\sqrt{ |G_{\alpha\epsilon_0}|}}{2(1+\sqrt{ |G_{\alpha\epsilon_0}|})}\bigg( \frac{7}{8}-\frac{1}{\sqrt{8\log  |G_{\alpha\epsilon_0}|}}\bigg)\ge C \alpha\epsilon_0,
\end{align*}
for some $C>0$ as long as $|G_{\alpha\epsilon_0}|\ge 2$. Hence,
\begin{align*}
\inf_{\widehat{\bU}}\sup_{\theta\in \mathcal{Y}(\CC,t,p_1,p_2,r)}\E_\theta d(\widehat{\bU}, \bU)&\gtrsim \inf_{\widehat{\bU}}\sup_{\theta\in \mathcal{X}'(t,\epsilon_0,\alpha\epsilon_0,\bU_0)}\E_\theta d(\widehat{\bU},\bU)\\
&\gtrsim  \bigg( \frac{\sigma\sqrt{4t^2+\sigma^2p_2}}{t^2}\sqrt{\log |G_{\alpha\epsilon_0}|}\land  \text{diam}(\CC)\bigg).
\end{align*}

\paragraph{Proof of Theorem \ref{lower.bnd.thm.2}.}
For some $\bU_0\in \CC$, similar to the proof of Theorem \ref{mdm.lower.thm}, we consider the $\delta$-packing set $G_\delta=G(B(\bU_0,\epsilon)\cap \CC,d,\delta)$,  where for any $\bU_i,\bU_j\in G_\delta$, $d(\bU_i,\bU_j)=\|\bU_i\bU_i^\top-\bU_j\bU_j^\top \|_F\ge \delta.$
Then, for given $t> 0$, we consider the subset $\mathcal{Z}'(t,\epsilon,\delta,\bU_0)=\{ (\bGam,\bU)\in\mathcal{Z}(\CC,t,p,r) : \bU \in G_\delta, \bGam =t{\bf I}_r\},$
so that $|\mathcal{Z}'(t,\epsilon,\delta,\bU_0)|=|G_\delta|.$
Let $P_i$ be the joint probability measure of $Y_k\sim_{i.i.d.} N(0,\bSig_i)$ with $k=1,...,n$ and $\bSig_i=t \bU_i\bU_i^\top+\sigma^2{\bf I}_p$. We have, for any $1\le i\ne j\le |G_\delta|$,
\begin{align*}
D(P_i,P_j) &=\frac{n}{2}\bigg( \text{tr}(\bSig_j^{-1}\bSig_i)-p+\log\bigg(\frac{\det \bSig_i}{\det \bSig_j}  \bigg) \bigg)\\
&=\frac{n}{2}\text{tr}\bigg( -\frac{t}{t+\sigma^2}\bU_i\bU_i^\top+\frac{t}{\sigma^2}\bU_j\bU_j^\top-\frac{t^2}{\sigma^2(t+\sigma^2)}\bU_i\bU_i^\top\bU_j\bU_j^\top \bigg)\\
&=\frac{nt^2}{2\sigma^2(\sigma^2+t)}(r-\|\bU_i^\top \bU_j\|_F^2) 
\le \frac{nt^2d^2(\bU_i,\bU_j)}{\sigma^2(\sigma^2+t)}
\le \frac{nt^2\epsilon^2}{\sigma^2(\sigma^2+t)},
\end{align*}
where the second equation follows from the Woodbury matrix identity and the second last inequality follows from Lemma \ref{dist.lem}. Now let $\epsilon=\epsilon_0$ and $\delta=\alpha\epsilon$ for some $\alpha\in(0,1)$. By assumption, 
\[
\bigg(\frac{c\sigma^2(\sigma^2+t)}{nt^2}\log |G_{\alpha\epsilon_0}| \land \text{diam}^2( \CC)\bigg)\le\epsilon_0^2\le \bigg(\frac{\sigma^2(\sigma^2+t)}{32nt^2}\log |G_{\alpha\epsilon_0}| \land \text{diam}^2( \CC)\bigg),
\]
for some $c\in(0,1/32)$.
It holds that $D(P_i,P_j) \le \frac{1}{16}\log |G_{\alpha\epsilon_0}|.$
Now by Lemma \ref{lower.lem}, it holds that, for $\theta=(\bGam,\bU)$,
\[
\inf_{\widehat{\bU}}\sup_{\theta\in \mathcal{Z}'(t,\epsilon_0,\alpha\epsilon_0,\bU_0)} P_\theta(d(\widehat{\bU},\bU )\ge \alpha\epsilon_0/2 )\ge \frac{\sqrt{  |G_{\alpha\epsilon_0}|}}{1+\sqrt{  |G_{\alpha\epsilon_0}|}}\bigg( \frac{7}{8}-\frac{1}{\sqrt{8\log   |G_{\alpha\epsilon_0}|}}\bigg).
\]
By Markov's inequality, as long as $ |G_{\alpha\epsilon_0}|\ge 2$, we have 
\begin{align*}
\inf_{\widehat{\bU}}\sup_{\theta\in \mathcal{Z}'(t,\epsilon_0,\alpha\epsilon_0,\bU)}\E_\theta d(\widehat{\bU}, \bU )&\ge C \alpha\epsilon_0,
\end{align*}
for some $C>0$. Therefore, since $\mathcal{Z}'(t,\epsilon_0,\alpha\epsilon_0,\bU)\subset \mathcal{Z}(\CC,t,p,r)$,
\begin{align*}
\inf_{\widehat{\bU}}\sup_{\theta\in \mathcal{Z}(\CC,t,p,r)}\mathcal{R}(\widehat{\bU}, \bU )&\ge \inf_{\widehat{\bU}}\sup_{\theta\in\mathcal{Z}'(t,\epsilon_0,\alpha\epsilon_0,\bU_0)}\E_\theta d(\widehat{\bU}, \bU )\\
&\gtrsim \bigg(\frac{{\sigma} \sqrt{\sigma^2+t}}{t\sqrt{n}}\sqrt{\log |G_{\alpha\epsilon_0}|}\land  \text{diam}( \CC)\bigg).
\end{align*}

\section{Calculation of Metric Entropies} \label{entropy.sec}

In this section, we prove the results in Section \ref{example.sec} by calculating metric entropies of some specific sets.
The calculation  relies on the following useful lemmas.

\bel[Varshamov-Gilbert Bound] \label{vg.bnd}
Let $\Omega=\{ 0,1 \}^n$ and $1\le d\le n/4$. Then there exists a subset $\{ \omega^{(1)},...,\omega^{(M)} \}$ of $\Omega$ such that $\|\omega^{(j)}\|_0=d$ for all $1\le j\le M$ and $\| \omega^{(j)}-\omega^{(k)}\|_0 \ge \frac{d}{2}$ for $0\le j<k\le M$, and $\log M\ge cd\log\frac{n}{d}$ where $c\ge 0.233.$
\eel 

The proof of the above version of Varshamov-Gilbert bound can be found, for example, in Lemma 4.10 in \cite{massart2007concentration}). The next two lemmas concern estimates of the covering/packing numbers of the orthogonal group.

\bel[\citealt{candes2011tight}] \label{plan.lem}
Define $\mathcal{P}_0=\{\bar{\bU}\bar{\bGam}\bar{\bV}^\top: \bar{\bU},\bar{\bV}\in O(p,2r), \|(\bar{\bGam}_{ii})_{1\le i\le 2r}\|_2=1\}$. Then for  any $\epsilon\in(0,\sqrt{2})$, there exists an $\epsilon$-covering set $H(\mathcal{P}_0,d_2,\epsilon)$ such that $| H(\mathcal{P}_0,d_2,\epsilon)| \le (c/\epsilon)^{2(2p+1)r}$ for some constant $c>0$.
\eel

\bel \label{entropy.g}
For any $V\in O(k,r)$, identifying the subspace $\text{span}(V)$ with its projection matrix $VV^\top$, define the metric on the Grassmannian manifold $G(k,r)$ by $\rho(VV^\top, UU^\top)=\|VV^\top-UU^\top\|_F$. Then for any $\epsilon\in (0,\sqrt{2(r\land (k-r))})$,
\[
\bigg( \frac{c_0}{\epsilon}\bigg)^{r(k-r)} \le \mathcal{N}(G(k,r),\rho,\epsilon)\le \bigg( \frac{c_1}{\epsilon}\bigg)^{r(k-r)},
\]
where $\mathcal{N}(E,\epsilon)$ is the $\epsilon$-covering number of $E$ and $c_0,c_1$ are absolute constants. Moreover, for any $V\in O(k,r)$ and any $\alpha\in(0,1)$, it holds that
\[
\bigg(\frac{c_0}{\alpha c_1} \bigg)^{r(k-r)}\le \mathcal{M}(\bB(V,\epsilon),\rho,\alpha\epsilon) \le  \bigg(\frac{2c_1}{\alpha c_0} \bigg)^{r(k-r)}.
\]
\eel

\begin{proof}
	We only prove the entropy upper bound
	\beq \label{51}
	\mathcal{M}(\bB(V,\epsilon),d,\alpha\epsilon) \le  \bigg(\frac{c_0}{\alpha c_1} \bigg)^{r(k-r)},
	\eeq
	as the other results has been proved in Lemma 1 of \cite{cai2013sparse}. Specifically, Let $G_{\epsilon}$ be the $\epsilon$-packing set of $O(k,r)$. It then holds that
	\begin{align*}
	&\mathcal{M}(O(k,r),d,\alpha\epsilon) \ge \sum_{V\in G_{\epsilon}} \mathcal{M}(\bB(V,\epsilon),d,\alpha\epsilon)\ge |G_\epsilon|\mathcal{M}(\bB(V^*,\epsilon),d,\alpha\epsilon)\\
	&=\mathcal{M}(O(k,r),d,\epsilon)\mathcal{M}(\bB(V^*,\epsilon),d,\alpha\epsilon)
	\end{align*}
	for some $V^*\in O(k,r)$. Hence,
	\[
	\mathcal{M}(\bB(V^*,\epsilon),d,\alpha\epsilon)\le \frac{	\mathcal{M}(O(k,r)),d,\alpha\epsilon)}{	\mathcal{M}(O(k,r)),d,\epsilon)}.
	\]
	By the equivalence between the packing and the covering numbers, it holds that
	\[
	\mathcal{M}(\bB(V^*,\epsilon),d,\alpha\epsilon)\le \frac{\mathcal{N}(O(k,r)),d,\alpha\epsilon/2)}{\mathcal{N}(O(k,r)),d,\epsilon)} \le \bigg(\frac{2c_1}{\alpha c_0} \bigg)^{r(k-r)},
	\]
	where the last inequality follows from the first statement of the lemma.
	Then (\ref{51}) holds since the metric $d$ is unitarily	invariant.
\end{proof}

The following lemma is an estimate of the Dudley's entropy integral for the orthogonal group $O(p,r)$.

\bel \label{plan.lem2}
For any given $\bU\in O(p,r)$, there exists some constant $C>0$ such that $\int_0^\infty \sqrt{\log  \mathcal{N}(\mathcal{T}(O(p,r),\bU),d_2,\epsilon)} d\epsilon \le C\sqrt{pr}$. Therefore, we have $\Delta^2(O(p,r))\le Cpr$.
\eel

\begin{proof}
	By definition, for any $\bG\in \mathcal{T}(O(p,r),\bU)$, it is at most rank $2r$, and suppose its SVD is $\bG=\bar{\bU}\bar{\bGam}\bar{\bV}^\top$, then $\bar{\bGam}$ is a diagonal matrix with nonnegative diagonal entries and Frobenius norm equal to one. Thus, if we define $\mathcal{P}_0=\{\bar{\bU}\bar{\bGam}\bar{\bV}^\top: \bar{\bU},\bar{\bV}\in O(p,2r), \|(\bar{\bGam}_{ii})_{1\le i\le 2r}\|_2=1\}$, then by Lemma \ref{szarek.lem},
	\[
	\mathcal{N}(\mathcal{T}(O(p,r),\bU),d_2,\epsilon) \le \mathcal{N}(\mathcal{P}_0,d_2,\epsilon).
	\]
	By Lemma \ref{plan.lem}, we can calculate that
	\begin{align}
	&\int_0^\infty \sqrt{\log  \mathcal{N}(\mathcal{T}(O(p,r),\bU),d_2,\epsilon)} d\epsilon\le \int_0^\infty \sqrt{\log  \mathcal{N}(\mathcal{P}_0,d_2,\epsilon)} d\epsilon\nonumber\\
	&\le C\sqrt{pr}\int_0^{\sqrt{2}} \sqrt{\log (c/\epsilon) }d\epsilon\le C\sqrt{pr}.
	\end{align}
	The second statement follows directly from the definition of $\Delta^2( O(p,r))$.
\end{proof}

\subsection{Sparse PCA/SVD: Proof of Proposition \ref{entropy.sparse} and Theorem \ref{sparse.thm}}

\paragraph{Matrix denoising model  with $\CC_S(p_1,r,k)$, or sparse SVD.} 
By Lemma \ref{vg.bnd}, we can construct a subset $\Theta_{\epsilon}(k)\subset \CC_S(p_1,r,k)$ as follows. Let $\Omega_M=\{\omega^{(1)},...,\omega^{(M)}\}\subset \{0,1\}^{p_1-r-1}$ be the set obtained from Lemma \ref{vg.bnd} where $n=p_1-r-1$, $d=k/e<(p_1-r-1)/4$ and $M$ is the smallest integer such that $\log M\ge cd\log n/d$, i.e., $M=\lceil \exp(ck\log \frac{e(p_1-r-1)}{k}) \rceil$. We define
\[
\Theta_\epsilon=\bigg\{  \begin{bmatrix}
\bv & \bold{0}\\
\bold{0} & {\bf I}_{r-1}
\end{bmatrix}:\bv= (\sqrt{1-\epsilon^2}, \epsilon \omega/\sqrt{d})\in \bS^{p_1-r-1}, \omega\in  \Omega_M\bigg\},\qquad \epsilon\in(0,1).
\]
Then $\Theta_\epsilon$ is a $\frac{\epsilon}{2}$-packing set of $\bB(\bU_0,\sqrt{2}\epsilon)\cap\CC_S(p_1,r,k)$ with $\bU_0=\begin{bmatrix}
\bv_0 & \bold{0}\\
\bold{0} & {\bf I}_{r-1}
\end{bmatrix}$ where $\bv_0=(1,0,...,0)^\top$, $|\Theta_\epsilon|= M$.
Now we set
\[
\epsilon^2=\frac{c_1(t^2+\sigma^2p_2)\sigma^2k\log (e(p_1-r-1)/k)}{t^4}\land 1,
\]
for some sufficiently small $c_1>0$.
It follows that
\[
\bigg(\frac{c_2\sigma^2(t^2+\sigma^2 p_2)}{t^4}\log |\Theta_\epsilon |\land 1\bigg)\le\epsilon^2\le \bigg(\frac{\sigma^2(t^2+\sigma^2p_2)}{640t^4}\log |\Theta_{\epsilon_0} |\land 1\bigg)
\]
for some $c_2\in(0,1/640)$.  So the condition of Theorem \ref{mdm.lower.thm} holds with $\epsilon_0=\sqrt{2}\epsilon$, $\alpha=1/(2\sqrt{2})$ and $\log |\Theta_\epsilon |\asymp k\log (ep_1/k)$.  Moreover, for any $\bU'\in O(k,r)$, suppose $M_\epsilon\subset O(k,r)$ is an $\alpha\epsilon$-packing set of $\bB(\bU',\epsilon)$ constructed as in Lemma \ref{entropy.g}, then the set
\beq \label{theta'}
\Theta'_{\epsilon}=\bigg\{ \bU=\begin{bmatrix}
	\bW\\
	\bold{0}
\end{bmatrix}, \bW\in M_\epsilon \bigg\} \subset \CC_S(p_1,r,k),
\eeq 
is an $\alpha\epsilon$-packing set of $\CC_S(p_1,r,k)\cap \bB(\bU_0,\epsilon)$ where $\bU_0=\begin{bmatrix} \bU'\\ \bold{0}\end{bmatrix}$, and $|\Theta'_{\epsilon}|\ge (c/\alpha)^{r(k-r)}$.  Now we set
\[
\epsilon^2=\frac{c_1(t^2+\sigma^2p_2)\sigma^2r(k-r)}{t^4}\land r^2,
\]
for some sufficiently small $c_1>0$.
It follows that
\[
\bigg(\frac{c_2\sigma^2(t^2+\sigma^2 p_2)}{t^4}\log |\Theta'_\epsilon |\land r\bigg)\le\epsilon^2\le \bigg(\frac{\sigma^2(t^2+\sigma^2p_2)}{640t^4}\log |\Theta'_{\epsilon_0} |\land r\bigg)
\]
for some $c_2\in(0,1/640)$. 
Thus, the condition of Theorem \ref{mdm.lower.thm} holds with $\log |\Theta'_\epsilon| \asymp r(k-r)$.

To obtain an upper bound for $\Delta(\CC_S(p_1,r,k))$, we notice that any element $\bH\in \CC_S(p_1,r,k)$ satisfies $\bH=\bH^\top$ and
\[
\max_{1\le i\le p_1}\|\bH_{i.}\|_0\le k,\quad\max_{1\le i\le p_1}\|\bH_{.i}\|_0\le k.
\]
Then $\mathcal{T}(\CC_S(p_1,r,k),\bU)$ can be covered by the union of its ${p_1\choose k}$ disjoint subsets, with each subset corresponding to a fixed sparsity configuration. Each of the above subsets can be identified with $\mathcal{T}(O(k,r),\bU')$ for some $\bU'\in O(k,r)$, and by Lemma \ref{plan.lem} and the proof of Lemma \ref{plan.lem2},  
\[
\mathcal{N}(\mathcal{T}(O(k,r),\bU'),d_2,\epsilon)\le (c/\epsilon)^{2r(2k+1)}.
\]
for any $\epsilon\in(0,\sqrt{2})$. Then by taking a union of the covering sets, we have
\[
\mathcal{N}(\mathcal{T}(\CC_S(p_1,r,k),\bU),d_2,\epsilon)\le {p_1\choose k}(c_1/\epsilon)^{2r(2k+1)}\le (ep_1/k)^k(c_1/\epsilon)^{2r(2k+1)}.
\]
As a result,
\begin{align*}
\int_0^\infty \sqrt{\log \mathcal{N}(\mathcal{T}(\CC_S(p_1,r,k),\bU),d_2,\epsilon)}d\epsilon&\le  \sqrt{2k\log(ep_1/k)}+\sqrt{2r(2k+1)}\int_0^{\sqrt{2}}\sqrt{\log\frac{c_1}{\epsilon}}d\epsilon\\
&\le C(\sqrt{k\log(ep_1/k)}+\sqrt{rk}).
\end{align*}
In addition, we also have
\[
\int_0^\infty \log \mathcal{N}(\mathcal{T}(\CC_S(p_1,r,k),\bU),d_2,\epsilon)d\epsilon\le C({k\log(ep_1/k)}+{rk}).
\]
The validity of Theorem \ref{mdm.risk.thm} reduces to the condition $\frac{t^2}{\sigma}\gtrsim {k\log(ep_1/k)}+{rk}.$ Note that when $r=O(1)$, this condition is satisfied whenever 
\[
\frac{\sigma\sqrt{t^2+\sigma^2 p_2}}{t^2} \bigg( \sqrt{k\log\frac{ep_1}{k}}+\sqrt{k} \bigg)\lesssim 1.
\]
In other words, in light of the minimax lower bound (from Theorem \ref{mdm.lower.thm}), whenever consistent estimation is possible, the condition $\frac{t^2}{\sigma}\gtrsim {k\log(ep_1/k)}+{k}$ is satisfied and the proposed estimator is minimax optimal. The final results follows by combining Theorems \ref{mdm.lower.thm} and \ref{mdm.risk.thm}.

\paragraph{Spiked Wishart model with $\CC_S(p,r,k)$, or sparse PCA.} We omitted the proof of this case as it is similar to the proof of the sparse SVD.

\subsection{Non-Negative PCA/SVD: Proof of Proposition \ref{entropy.NN} and Theorem \ref{NN.thm}}

\paragraph{Matrix denoising model with $\CC_{N}(p_1,r)$, or non-negative SVD.} On the one hand, with Lemma \ref{vg.bnd}, we can construct a subset $\Theta_{\epsilon}\subset O(p_1,r)$ as follows. Let $\Omega_M=\{\omega^{(1)},...,\omega^{(M)}\}\subset \{0,1\}^{n}$ be the set obtained from Lemma \ref{vg.bnd} where $n=p_1-r-1$, $d=(p_1-r-1)/4$ and $M$ is the smallest integer such that $\log M\ge cd\log n/d$, i.e., $M=\lceil \exp(\frac{c(p_1-r-1)\log 2}{2}) \rceil$. Following the idea of \cite{vu2012minimax} and \cite{cai2013sparse}, we define
\[
\Theta_\epsilon=\bigg\{  \begin{bmatrix}
\bv & \bold{0}\\
\bold{0} & {\bf I}_{r-1}
\end{bmatrix}:\bv= (\sqrt{1-\epsilon^2}, \epsilon \omega/\sqrt{d})\in \bS^{p_1-r-1}, \omega\in  \Omega_M\bigg\},\qquad \epsilon\in(0,1).
\]
Then it holds that $\Theta_\epsilon\subset \bB(\bU_0,\sqrt{2}\epsilon)$ for $\bU_0=\begin{bmatrix}
\bv_0 & \bold{0}\\
\bold{0} & {\bf I}_{r-1}
\end{bmatrix}$ where $\bv_0=(1,0,...,0)^\top$, $|\Theta_\epsilon|= M$, and that for any $\bU\ne \bU'\in\Theta_\epsilon$, 
\[
d(\bU,\bU')\ge\sqrt{2}\cdot\sqrt{1-(1-\epsilon^2/8)^2}\ge \frac{\epsilon}{2}.
\]
In other words, $\Theta_\epsilon$ is a $\frac{\epsilon}{2}$-packing set of $\bB(\bU_0,\sqrt{2}\epsilon)\cap\CC_{NN}(p_1,r)$.
Now we set
\[
\epsilon^2=\frac{c_1(t^2+\sigma^2p_2)\sigma^2(p_1-r-1)}{t^4}\land 1,
\]
for some sufficiently small $c_1>0$.
It follows that
\[
\bigg(\frac{c_2\sigma^2(t^2+\sigma^2p_2)}{t^4}\log |\Theta_\epsilon |\land 1\bigg)\le\epsilon^2\le \bigg(\frac{\sigma^2(t^2+\sigma^2p_2)}{640t^4}\log |\Theta_{\epsilon_0} |\land 1\bigg)
\]
for some $c_2\in(0,1/640)$.  So the condition of Theorem \ref{mdm.lower.thm} holds with $\epsilon_0=\sqrt{2}\epsilon$, $\alpha=1/(2\sqrt{2})$ and $\log |\Theta_\epsilon |\asymp p_1$.

On the other hand, we need to obtain an upper bound for $\Delta(\CC_{N}(p_1,r))$.
To bound the Dudley's entropy integral $\int_0^\infty \sqrt{\log  \mathcal{N}(\mathcal{T}(\CC_{N}(p_1,r),\bU),d_2,\epsilon)} d\epsilon$, we simply use the fact that $\CC_{N}(p_1,r)\subset O(p_1,r)$ and 
\[
\mathcal{N}(\mathcal{T}(\CC_{N}(p_1,r),\bU),d_2,\epsilon) \le \mathcal{N}(\mathcal{T}(O(p_1,r),\bU), d_2,\epsilon).
\]
Then by Lemma \ref{plan.lem2}, we have $\Delta^2(\CC_{NN}(p_1,r))\lesssim {p_1r}$. Combining Theorems \ref{mdm.lower.thm} and \ref{mdm.risk.thm}, we have $\Delta^2(\CC_{NN}(p_1,r))\gtrsim \log |\Theta_\epsilon|$, which implies $\Delta^2(\CC_{NN}(p_1,r))\asymp\log |\Theta_\epsilon|\asymp p_1$ if $r=O(1)$. Again, Theorem 8 requires $\frac{t^2}{\sigma}\gtrsim rp_1.$ Note that when $r=O(1)$, this condition is satisfied whenever 
\[
\frac{\sigma\sqrt{p_1(t^2+\sigma^2 p_2)}}{t^2}\lesssim 1.
\]
In other words, in light of the minimax lower bound (from Theorem \ref{mdm.lower.thm}), whenever consistent estimation is possible, the condition $\frac{t^2}{\sigma}\gtrsim p_1$ is satisfied and the proposed estimator is minimax optimal. 

\paragraph{Spiked Wishart model with $\CC_N(p,r)$, or non-negative PCA.} Similarly, let $\Omega_M=\{\omega^{(1)},...,\omega^{(M)}\}\subset \{0,1\}^{p-r-1}$ be the set obtained from Lemma \ref{vg.bnd} where $d=(p-r-1)/4$ and $M$ is the smallest integer such that $\log M\ge cd\log (p-r-1)/d$, i.e., $M=\lceil \exp(\frac{c(p-r-1)\log 2}{2}) \rceil$. We define
\[
\Theta_\epsilon=\bigg\{  \begin{bmatrix}
\bv & \bold{0}\\
\bold{0} & {\bf I}_{r-1}
\end{bmatrix}:\bv= (\sqrt{1-\epsilon^2}, \epsilon \omega/\sqrt{d})\in \bS^{p-r-1}, \omega\in  \Omega_M\bigg\},\qquad \epsilon\in(0,1).
\]
Then it holds that $\Theta_\epsilon\subset \bB(\bU_0,\sqrt{2}\epsilon)$ for $\bU_0=\begin{bmatrix}
\bv_0 & \bold{0}\\
\bold{0} & {\bf I}_{r-1}
\end{bmatrix}$ where $\bv_0=(1,0,...,0)^\top$, $|\Theta_\epsilon|= M$, and that for any $\bU\ne \bU'\in\Theta_\epsilon$, 
\[
d(\bU,\bU')\ge\sqrt{2}\cdot\sqrt{1-(1-\epsilon^2/8)^2}\ge \frac{\epsilon}{2}.
\]
In other words, $\Theta_\epsilon$ is a $\frac{\epsilon}{2}$-packing set of $\bB(\bU_0,\sqrt{2}\epsilon)\cap\CC_{NN}(p,r)$.
Now we set
\[
\epsilon^2=\frac{c_1\sigma^2(\sigma^2+t)(p-r-1)}{nt^2}\land 1,
\]
for some sufficiently small $c_1>0$.
It follows that
\[
\bigg(\frac{c_2\sigma^2(\sigma^2+t)}{nt^2}\log |\Theta_\epsilon |\land 1\bigg)\le\epsilon^2\le \bigg(\frac{\sigma^2(\sigma^2+t)}{nt^2} \frac{(p-r-1)\log 2}{10}\land 1\bigg) \le \bigg(\frac{\sigma^2(\sigma^2+t)}{32nt^2}\log |\Theta_\epsilon |\land 1\bigg)
\]
for some $c_2\in(0,1/32)$, so that condition of Theorem \ref{lower.bnd.thm.2} holds and $\log |\Theta_\epsilon |\asymp p$. The rest of the arguments such as the calculation of Dudley's entropy integral are the same as the above proof of the non-negative SVD.

\subsection{Subspace PCA/SVD: Proof of Proposition \ref{entropy.sub} and Theorem \ref{sub.thm}}

To prove this proposition, in light of Lemmas \ref{plan.lem}, \ref{entropy.g} and \ref{plan.lem2}, it suffices to establish the isometry between $(\CC_A(p,r,k),d)$ and $(O(k,r), d)$. Let $\bQ\in O(p,k)$ has its columns being the basis of the null space of $A$. We consider the map $F: O(k,r)\to\CC_A(p,r,k)$ where $F(\bW)=\bQ\bW$. To show that $F$ is a bijection, we notice that 
\begin{enumerate}
	\item For any $\bG\in \CC_A(p,r,k)$, for each of its columns $\bQ_{.i}$, there exists some $\bv_i\in \bS^{k-1}$ such that $\bG_{.i}=\bQ\bv_i$ and $\bv_i^\top\bv_j=\bv_i^\top \bQ^\top\bQ\bv_j=\bG_{.i}^\top\bG_{.j}=0$. Then let $\bW=[\bv_1,...,\bv_r]\in O(k,r)$, apparently, we have $F(\bW)=\bG$. This proves that the map is onto.
	\item For any $\bW_1\ne \bW_2\in O(k,r)$, it follows that $F(\bW_1)\ne F(\bW_2)$. This proves the injection.
\end{enumerate}
To show the map $F$ is isometric, we notice that
\begin{enumerate}
	\item For any $\bG_1=F(\bW_1),\bG_2=F(\bW_2)\in \CC_A(p,r,k)$,
	\begin{align*}
	d(F(\bW_1),F(\bW_2))&=\|\bQ\bW_1\bW_1^\top \bQ^\top-\bQ\bW_2\bW_2^\top \bQ^\top\|_F\\
	&\le \|\bQ\|^2\|\bW_1\bW_1^\top-\bW_2\bW_2^\top\|_F\\
	&\le d(\bW_1,\bW_2).
	\end{align*}
	\item For any $\bW_1,\bW_2\in O(k,r)$, 
	\[
	d(\bW_1,\bW_2)=\|\bQ^\top\bQ\bW_1\bW_1^\top \bQ^\top-\bQ\bW_2\bW_2^\top \bQ^\top\bQ\|_F\le 	d(F(\bW_1,\bW_2)).
	\]
\end{enumerate}
Thus $d(F(\bW_1),F(\bW_2))=d(\bW_1,\bW_2)$. 

\subsection{Spectral Clustering: Proof of Proposition \ref{entropy.pm} and Theorem \ref{pm.thm}}

The upper bound $\Delta^2(\CC_{\pm}^n)\lesssim n$ follows from the same argument as in the proof of Proposition 15. For the second statement, by Lemma \ref{vg.bnd}, we can construct a subset $\Theta(d)\subset \bS^{n-1}$ as follows. Let $\Omega_M=\{\omega^{(1)},...,\omega^{(M)}\}\subset \{0,1\}^{n}$ be the set obtained from Lemma \ref{vg.bnd} where $\|\omega^{(j)}\|_0=d\le n/4$ for all $1\le j\le n$ and $M$ is the smallest integer such that $\log M\ge cd$, i.e., $M=\lceil \exp(cd\log \frac{n}{d}) \rceil$. We define
\[
\Theta(d)=\bigg\{  \frac{2|\omega-0.5\cdot{\bf 1}|}{\sqrt{n}}\in \CC_{\pm}^n: \omega\in  \Omega_M\cup \{ (0,...,0) \}\bigg\},
\]
where ${\bf 1}=(1,...,1)^\top\in\R^n$.
Then since for $\bu_0=(-1/\sqrt{n},...,-1/\sqrt{n})^\top$ and any $\bu\in \Theta(d)$,
\[
d(\bu_0,\bu)\le \|\bu_0-\bu\|_2\le 2\sqrt{\frac{d}{n}},
\] 
it holds that $\Theta(d)\subset \bB(\bu_0, 2\sqrt{d/n})$ with and that for any $\bu\ne \bu'\in \Theta(d)$,
\[
d(\bu,\bu')\ge \frac{1}{\sqrt{2}}\|\bu-\bu'\|_2\ge \sqrt{\frac{d}{n}}
\]
so that $\Theta(d)$ is a $\sqrt{\frac{d}{n}}$-packing set of $\bB(\bu_0,2\sqrt{d/n})\cap\CC_{\pm}^{n}$.
Now since $t^2=C\sigma^2(n+\sqrt{np})$, we can set
\[
\epsilon_0=\sqrt{\frac{d}{n}},\quad \text{where} \quad d=c_1n,
\]
for some sufficiently small $c_1>0$, and thus it follows that
\[
\bigg(\frac{c_2\sigma^2(t^2+\sigma^2p)}{t^4}\log |\Theta(d) |\land 1\bigg)\le\epsilon_0^2\le \bigg(\frac{\sigma^2(t^2+\sigma^2p)}{128t^4}\log |\Theta(d) |\land 1\bigg)
\]
for some $c_2\in(0,1/128)$. 
So the condition of Theorem 5 holds with $\alpha=1/2$ and $\log |\Theta(d) |\asymp n$.

\section{Proof of Technical Lemmas} \label{lem.sec}

\paragraph{Proof of Lemma \ref{F.inn.lem}.} The first inequality can be proved by
\begin{align*}
\langle \bU\bGam^2\bU, \bU\bU^\top-\bW\bW^\top\rangle &=\text{tr}(\bU\bGam^2\bU^\top)-\text{tr}(\bW^\top \bU\bGam^2\bU^\top \bW)\\
&=\text{tr}(\bGam^2)-\text{tr}(\bGam^2\bU^\top \bW\bW^\top \bU)\\
&=\sum_{i=1}^r \lambda_i^2(1-(\bU^\top\bW\bW^\top\bU)_{ii})\\
&\ge \lambda_r^2(r-\text{tr}(\bU^\top\bW\bW^\top\bU))\\
&=\frac{\lambda_r^2}{2}\|\bU\bU^\top-\bW\bW^\top\|_F^2.
\end{align*}
The other inequality follows from the same rationale.

\paragraph{Proof of Lemma \ref{entropy.chaos.lem}.}
Throughout the proof, for simplicity, we write $\cP=\cP(\CC,\bU)$ and $\mathcal{T}=\mathcal{T}(\CC,\bU)$.
By Corollary 2.3.2 of \cite{talagrand2014upper}, for any metric space $(T,d)$, if we define
\beq \label{en}
e_n(T)=\inf\{ \epsilon: \mathcal{N}(T,d,\epsilon)\le N_n\},\quad \text{where $N_0=1; N_n=2^{2^n}$ for $n\ge 1,$}
\eeq
then there exists some constant $K(\alpha)$ only depending on $\alpha$ such that
\beq \label{gamma-en}
\gamma_\alpha(T,d)\le K(\alpha)\sum_{n\ge 0} 2^{n/\alpha}e_n(T).
\eeq
The following inequalities establish the correspondence between $e_n$ and the Dudley's entropy integral, 
\beq \label{en-dudley}
\begin{aligned}
	&\sum_{n\ge 0} 2^{n/2}e_n(T)\le C\int_0^{\infty} \sqrt{\log \mathcal{N}(T,d,\epsilon)}d\epsilon,\\
	&\sum_{n\ge 0} 2^{n}e_n(T)\le C\int_0^{\infty} {\log \mathcal{N}(T,d,\epsilon)}d\epsilon,
\end{aligned}
\eeq
whose derivation is delayed to the end of this proof.
Combining (\ref{gamma-en}) and (\ref{en-dudley}), it follows that
\beq  \label{gamma-dudley}
\gamma_\alpha(T,d)\le K(\alpha)\int_0^{\infty} {\log^{1/\alpha} \mathcal{N}(T,d,\epsilon)}d\epsilon.
\eeq
By (\ref{gamma-dudley}), it suffices to obtain estimates of the metric entropies ${\log \mathcal{N}(\cP,d_\infty,\epsilon)}$ and $\sqrt{\log \mathcal{N}(\cP,d_2,\epsilon)}$. By definition of $\mathcal{T}$, apparently $(\cP,d_\infty)$ is isomorphic to $(\mathcal{T},d_\infty)$,
then by Lemma \ref{szarek.lem}, it holds that
\[
\mathcal{N}(\cP,d_\infty,\epsilon)= \mathcal{N}(\mathcal{T},d_\infty,\epsilon).
\]
Along with the fact that, for any $\bG_1,\bG_2\in \mathcal{T}$, $d_\infty(\bG_1,\bG_2)\le d_2(\bG_1,\bG_2)$ and therefore
\[
\mathcal{N}(\mathcal{T},d_\infty,\epsilon)\le  \mathcal{N}(\mathcal{T},d_2,\epsilon),
\]
we prove the first statement of the lemma.
On the other hand, consider the map $F: (\cP,d_2)\to (\mathcal{T},d_2)$ where for any $\bD\in \cP$, $F(\bD)\in \R^{p_1\times p_1}$ is the submatrix of $\bD$  by extracting its entries in the first $p_1$ columns and rows. Then, for any $\bD_1,\bD_2\in \cP$, it holds that
\[
d_2(F(\bD_1),F(\bD_2))=\| F(\bD_1)-F(\bD_2)\|_F=\frac{1}{\sqrt{p_2}} d_2(\bD_1,\bD_2).
\] 
Again, applying Lemma 6, we have
\[
\mathcal{N}(\cP,d_2,\epsilon)= \mathcal{N}(\mathcal{T},d_2,\epsilon/\sqrt{p_2}).
\]
The second statement of the lemma then follows simply from the change of variable
\[
\gamma_2(\cP,d_2)\le C_2\int_0^{\infty}\sqrt{\log\mathcal{N}(\mathcal{T},d_2,\epsilon/\sqrt{p_2})}d\epsilon= C_2\sqrt{p_2}\int_0^{\infty}\sqrt{\log\mathcal{N}(\mathcal{T},d_2,\epsilon)}d\epsilon.
\]
\paragraph{Proof of (\ref{en-dudley}).} The proof of the first inequality can be found, for example, on page 22 of \cite{talagrand2014upper}. Nevertheless, we provide a detailed proof for completeness. By definition of $e_n$, if $\epsilon<e_n(T)$, we have $\mathcal{N}(T,d,\epsilon)>N_n$ and $\mathcal{N}(T,d,\epsilon)\ge N_n+1$. Then
\[
\sqrt{\log(1+N_n)}(e_n(T)-e_{n+1}(T))\le \int_{e_{n+1}(T)}^{e_n(T)} \sqrt{\log \mathcal{N}(T,d,\epsilon)}.
\]
Since $\log (1+N_n)\ge 2^n\log 2$ for $n\ge 0$, summation over $n\ge 0$ yields
\[
\sqrt{\log 2} \sum_{n\ge 0}2^{n/2}(e_n-e_{n+1}(T))\le \int_0^{e_0(T)} \sqrt{\log \mathcal{N}(T,d,\epsilon)}.
\]
Then the final inequality (\ref{en-dudley}) follows by noting that
\begin{align*}
\sum_{n\ge 0}2^{n/2}(e_n-e_{n+1}(T))&= \sum_{n\ge 0}2^{n/2}e_n(T)-\sum_{n\ge 1}2^{(n-1)/2}e_n(T)\\
&\ge (1-1/\sqrt{2}) \sum_{n\ge 0}2^{n/2}e_n(T).
\end{align*}
The second inequality can be obtained similarly by working with the inequality
\[
{\log(1+N_n)}(e_n(T)-e_{n+1}(T))\le \int_{e_{n+1}(T)}^{e_n(T)} {\log \mathcal{N}(T,d,\epsilon)}.
\]

\paragraph{Proof of Lemma \ref{kl.mix.prop}.} 
The proof of this lemma generalizes the ideas in \cite{cai2018rate} and \cite{ma2019optimalb}. In general, direct calculation of $D({P}_i,P_j)$ is difficult. We detour by introducing an approximate density of $P_i$ as
\begin{align*}
\tilde{P}_i(\bY)=\frac{\sigma^{-p_1p_2}}{(2\pi)^{p_1p_2/2}}\int&\exp(-\|\bY-t\bU_i\bW^\top\|_F^2/(2\sigma^2))\bigg( \frac{p_2}{2\pi}\bigg)^{rp_2/2}\exp(-p_2\|\bW\|_F^2/2)d\bW.
\end{align*}
Now for $\bY\sim \tilde{P}_i$, if $Y_k$ is the $k$-th column of $\bY$, we have
\beq \label{Y_i.dist}
Y_k| \bU_i \sim_{i.i.d.} N\bigg( 0, \sigma^2\bigg(I_n-\frac{4t^2}{4t^2+\sigma^2p_2}\bU_i \bU_i^\top\bigg)^{-1} \bigg) = N\bigg(0,\sigma^2I_n+\frac{4t^2}{p_2}\bU_i\bU_i^\top\bigg), 
\eeq
for $k=1,...,p_2.$ 
It is well-known that the KL-divergence between two $p$-dimensional multivariate Gaussian distribution is
\[
D( N(\mu_0,\bSig_0)\| N(\mu_1,\bSig_1)) =\frac{1}{2}\bigg(  \text{tr}(\bSig_0^{-1}\bSig_1)+(\mu_1-\mu_0)^\top \bSig_1^{-1}(\mu_1-\mu_0)-p+\log\bigg( \frac{\det \bSig_1}{\det \bSig_0}\bigg) \bigg).
\]
As a result, we can calculate that for any $\tilde{P}_i$ and $\tilde{P}_j$,
\begin{align} \label{KL.tilde}
D(\tilde{P}_i, \tilde{P}_j) &=\frac{p_2}{2}\bigg\{ \text{tr}\bigg( \bigg( I_{p_1} -\frac{4t^2}{4t^2+\sigma^2p_2}\bU_i\bU_i^\top \bigg)\bigg( I_{p_1} +\frac{4t^2}{\sigma^2p_2}\bU_j \bU_j^\top \bigg)\bigg)-p_1 \bigg\}\nonumber \\
&\le \frac{Ct^4}{4t^2+\sigma^2p_2}(r-\|\bU_i^\top \bU_j\|_F^2)\nonumber\\
&=\frac{Ct^4d(\bU_i,\bU_j)}{4t^2+\sigma^2p_2}
\end{align}
where the last inequality follows from Lemma \ref{dist.lem}.
Hence, the proof of this proposition is complete if we can show that there exist some constant $C>0$ such that
\beq \label{affinity.2}
D({P}_i,{P}_j)\le D(\tilde{P}_i,\tilde{P}_j)+C.
\eeq
The rest of the proof is devoted to the proof of (\ref{affinity.2}).
\paragraph{Proof of (\ref{affinity.2}).} Define the event $\mathcal{G}=\{\bW\in \R^{r\times p_2}: 1/2\le \lambda_{\min}(\bW)\le \lambda_{\max}(\bW)\le 2\}$. For any given $u$,
\begin{align} \label{pu/qu}
\frac{{P}_i}{\tilde{P}_i}&=\frac{1}{(2\pi)^{\frac{rp_2}{2}}(\frac{\sigma^2}{4t^2+\sigma^2p_2})^{\frac{rp_2}{2}}}\exp\bigg(\frac{1}{2\sigma^2} \sum_{k=1}^{p_2}Y_{k}^\top(I_{p_1}-\frac{4t^2}{4t^2+\sigma^2p_2}\bU_i\bU_i^\top) Y_{k}\bigg) \nonumber \\
&\quad\times C_{\bU_i,t} \int_{\mathcal{G}}\exp(-\|\bY-t\bU_i\bW^\top\|_F^2/(2\sigma^2)-p_2\|\bW\|_F^2/2)d\bW \nonumber \\
&=\bigg(\frac{4t^2+\sigma^2p_2}{2\pi\sigma^2}\bigg)^{p_2r/2}\exp\bigg(-(4t^2+\sigma^2p_2)\bigg\|\bW-\frac{2t}{4t^2+\sigma^2p_2}\bU_i^\top\bY\bigg\|_F^2/2\bigg)d\bW\nonumber \\
&=C_{\bU_i,t}P\bigg(\bW'\in\mathcal{G}\bigg|  \bW' \sim N\bigg( \frac{2t}{4t^2+\sigma^2p_2}\bU_i^\top\bY, \frac{\sigma^2}{4t^2+\sigma^2p_2}{\bf I}_{p_1} \bigg)\bigg)\nonumber \\
&\le C_{\bU_i,t}.
\end{align}
Recall that
\[
C_{\bU_i,t}^{-1}= P\big( \bW=(w_{jk})\in\mathcal{G}| w_{jk}\sim N(0,1/p_2)  \big).
\]
By concentration of measure inequalities for Gaussian random matrices (see, for example, Corollary 5.35 of \cite{vershynin2010introduction}), we have, for sufficiently large $(p_2,r)$,
\beq \label{key.ineq}
P(\bW\in\mathcal{G})\ge 1-2\exp(-cp_2),
\eeq
for some constant $c>0$.
In other words, we have
\beq 
C^{-1}_{\bU_i,t}\ge 1-p_2^{-c}
\eeq
and 
\beq
\frac{{P}_i}{\tilde{P}_i}\le 1+p_2^{-c}
\eeq
uniformly for some constant $c>0.$ Thus, for some constant $\delta>0$, we have
\begin{align}
D({P}_i,{P}_j)&=\int {P}_i\bigg[\log \bigg( \frac{{P}_i}{\tilde{P}_i} \bigg)+\log \bigg( \frac{\tilde{P}_i}{\tilde{P}_j} \bigg)+\log \bigg( \frac{\tilde{P}_j}{{P}_j} \bigg)\bigg]d\bY\nonumber\\
&\le \log(1+\delta)+D(\tilde{P}_i, \tilde{P}_j) +\int ({P}_i-\tilde{P}_i)\log \bigg( \frac{\tilde{P}_i}{\tilde{P}_j} \bigg)d\bY+\int{P}_i\log \bigg( \frac{\tilde{P}_i}{{P}_j} \bigg)d\bY\nonumber \\
&\le \log(1+\delta)+D(\tilde{P}_i, \tilde{P}_j)  +\int\tilde{P}_i \bigg(\frac{{P}_i}{\tilde{P}_i}-1\bigg)\log \bigg( \frac{\tilde{P}_i}{\tilde{P}_j} \bigg)d\bY\nonumber\\
&\quad+(1+\delta)\int\tilde{P}_i\bigg|\log \bigg( \frac{\tilde{P}_j}{{P}_j} \bigg)\bigg|d\bY\nonumber\\
&\le \log(1+\delta)+D(\tilde{P}_i, \tilde{P}_j)  +p_2^{-c}\int\tilde{P}_i \bigg|\log \bigg( \frac{\tilde{P}_i}{\tilde{P}_j} \bigg)\bigg|d\bY+(1+\delta)\int\tilde{P}_i\bigg|\log \bigg( \frac{\tilde{P}_j}{{P}_j} \bigg)\bigg|d\bY. \label{eq}
\end{align}
Now since
\begin{align*}
\int\tilde{P}_i \bigg|\log \bigg( \frac{\tilde{P}_i}{\tilde{P}_j} \bigg)\bigg|d\bY&=\frac{1}{2\sigma^2}\int\tilde{P}_i \bigg|\frac{4t^2}{4t^2+\sigma^2p_2}\sum_{k=1}^{p_2}Y_{k}^\top(\bU_i\bU_i^\top-\bU_j\bU_j^\top)Y_{k}\bigg|d\bY\\
&\le \frac{1}{2\sigma^2}\E\bigg[ \frac{4t^2}{4t^2+\sigma^2p_2}\sum_{k=1}^{p_2}Y_{k}^\top(\bU_i\bU_i^\top+\bU_j\bU_j^\top)Y_{k} \bigg]\\
&=\frac{4t^2p_2}{2\sigma^2(4t^2+\sigma^2p_2)}\text{tr}\bigg( (\bU_i\bU_i^\top+\bU_j\bU_j^\top)\big(\sigma^2I_{p_1}+\frac{4t^2}{p_2}\bU_i\bU_i^\top\big)  \bigg)\\
&\le \frac{4t^2p_2}{4t^2+\sigma^2p_2}\text{tr}\bigg(\bU_i^\top\big(I_{p_1}+\frac{4t^2}{\sigma^2p_2} \big)\bU_i  \bigg)\\
&=\frac{4rt^2}{\sigma^2}\le  rp_2,
\end{align*}
where in the second row the expectation is with respect to $Y_{k}\sim N\big(0,\sigma^2I_{p_1}+\frac{4t^2\sigma^2}{\sigma^2p_2}\bU_i\bU_i^\top \big)$.
we know that the third term in (\ref{eq}) can be bounded by
\[
p_2^{-c}\int\tilde{P}_i \bigg|\log \bigg( \frac{\tilde{P}_i}{\tilde{P}_j} \bigg)\bigg|d\bY\le rp_2\cdot p_2^{-c}\le C
\]
for some constants $C,c>0$.
Finally, by (\ref{pu/qu}), we have
\begin{align*}
\int \tilde{P}_i\bigg|\log \bigg( \frac{\tilde{P}_j}{{P}_j} \bigg)\bigg|dY&\le \int\tilde{P}_i\bigg|\log \frac{1}{C_{\bU_j,t}}\bigg|dY+\int \tilde{P}_i\bigg|\log \frac{1}{P(\bW'\in\mathcal{G}|E)}\bigg|d\bY,
\end{align*}
where we denoted 
\[
E=\bigg\{ \bW' \sim N\bigg( \frac{2t}{4t^2+\sigma^2p_2}\bU_i^\top\bY, \frac{\sigma^2}{4t^2+\sigma^2p_2}{\bf I}_{p_1} \bigg) \bigg\}.
\]
Now on the one hand,
\[
\int\tilde{P}_i\bigg|\log \frac{1}{C_{\bU_i,t}}\bigg|d\bY\le \big(\log(1+\delta) \lor |\log(1-\delta)^{-1}|\big).
\]
On the other hand, for fixed $\bY$ and $\bU_i^\top\bY\in\R^{r\times p_2}$, we can find $\bQ\in O(p_2,p_2-r)$ which is orthogonal to $\bU_i^\top\bY$, i.e., $\bU_i^\top\bY\bQ=0$. Then $\bW'\bQ\in \R^{r\times (p_2-r)}$ are i.i.d. normal distributed with mean $0$ and variance $\frac{\sigma^2}{4t^2+\sigma^2p_2}$. Then again by standard result in random matrix (e.g. Corollary 5.35 in \cite{vershynin2010introduction}), we have
\[
\lambda_{\min}(\bW')=\lambda_r(\bW')\ge \lambda_r(\bW'\bQ)\ge \frac{\sigma}{\sqrt{4t^2+\sigma^2p_2}}(\sqrt{p_2-r}-\sqrt{r}-x)
\] 
with probability at least $1-2\exp(-x^2/2)$. Since $t^2<\sigma^2 p_2/4$, for $p_2$ sufficiently large, we can find $c$ such that by setting $x=c\sqrt{p_2}$, 
\beq \label{G1}
P( \lambda_{\min}(\bW')\ge 1/2) \ge 1-e^{-cp_2}.
\eeq 
Analogous to the argument on $\lambda_{\min}(\bW')$, we also have
\beq
P( \lambda_{\max}(\bW')\le 2) \ge 1-e^{-cp_2}.
\eeq
Thus, by the union bound inequality, we have
\[
P(\bW'\in\mathcal{G}) \ge 1-2e^{-cp_2},
\]
and consequently,
\[
\int \tilde{P}_i\bigg|\log \frac{1}{P(\bW'\in\mathcal{G}|E)}\bigg|d\bY\le \bigg|\log \frac{1}{1-p_2^{-c}}\bigg|\le p_2^{-c}.
\]
This helps us to bound the last term of (\ref{eq}).
Combining the above results, we have proven the inequality (\ref{affinity.2}) and therefore completed the proof.

\vskip 0.2in
\bibliography{reference}

\begin{thebibliography}{75}
\providecommand{\natexlab}[1]{#1}
\providecommand{\url}[1]{\texttt{#1}}
\expandafter\ifx\csname urlstyle\endcsname\relax
  \providecommand{\doi}[1]{doi: #1}\else
  \providecommand{\doi}{doi: \begingroup \urlstyle{rm}\Url}\fi

\bibitem[Arcones and Gin{\'e}(1993)]{arcones1993decoupling}
Miguel~A Arcones and Evarist Gin{\'e}.
\newblock On decoupling, series expansions, and tail behavior of chaos
  processes.
\newblock \emph{J. Theor. Probab.}, 6\penalty0 (1):\penalty0 101--122, 1993.

\bibitem[Azizyan et~al.(2013)Azizyan, Singh, and Wasserman]{azizyan2013minimax}
Martin Azizyan, Aarti Singh, and Larry Wasserman.
\newblock Minimax theory for high-dimensional gaussian mixtures with sparse
  mean separation.
\newblock In \emph{NIPS}, pages 2139--2147, 2013.

\bibitem[Bai and Yao(2008)]{bai2008central}
Zhidong Bai and Jian-feng Yao.
\newblock Central limit theorems for eigenvalues in a spiked population model.
\newblock In \emph{Annales de l'IHP Probabilit{\'e}s et Statistiques},
  volume~44, pages 447--474, 2008.

\bibitem[Baik and Silverstein(2006)]{baik2006eigenvalues}
Jinho Baik and Jack~W Silverstein.
\newblock Eigenvalues of large sample covariance matrices of spiked population
  models.
\newblock \emph{J. Multiv. Anal.}, 97\penalty0 (6):\penalty0 1382--1408, 2006.

\bibitem[Bandeira et~al.(2017)Bandeira, Boumal, and
  Singer]{bandeira2017tightness}
Afonso~S Bandeira, Nicolas Boumal, and Amit Singer.
\newblock Tightness of the maximum likelihood semidefinite relaxation for
  angular synchronization.
\newblock \emph{Mathematical Programming}, 163\penalty0 (1-2):\penalty0
  145--167, 2017.

\bibitem[Bao et~al.(2018)Bao, Ding, and Wang]{bao2018singular}
Zhigang Bao, Xiucai Ding, and Ke~Wang.
\newblock Singular vector and singular subspace distribution for the matrix
  denoising model.
\newblock \emph{arXiv preprint arXiv:1809.10476}, 2018.

\bibitem[Bartlett and Mendelson(2002)]{bartlett2002rademacher}
Peter~L Bartlett and Shahar Mendelson.
\newblock Rademacher and gaussian complexities: Risk bounds and structural
  results.
\newblock \emph{J. Mach. Learn. Res.}, 3\penalty0 (Nov):\penalty0 463--482,
  2002.

\bibitem[Birnbaum et~al.(2013)Birnbaum, Johnstone, Nadler, and
  Paul]{birnbaum2013minimax}
Aharon Birnbaum, Iain~M Johnstone, Boaz Nadler, and Debashis Paul.
\newblock Minimax bounds for sparse pca with noisy high-dimensional data.
\newblock \emph{Ann. Statist.}, 41\penalty0 (3):\penalty0 1055--1084, 2013.

\bibitem[Boumal(2016)]{boumal2016nonconvex}
Nicolas Boumal.
\newblock Nonconvex phase synchronization.
\newblock \emph{SIAM J. Optimiz.}, 26\penalty0 (4):\penalty0 2355--2377, 2016.

\bibitem[Bousquet et~al.(2002)Bousquet, Koltchinskii, and
  Panchenko]{bousquet2002some}
Olivier Bousquet, Vladimir Koltchinskii, and Dmitriy Panchenko.
\newblock Some local measures of complexity of convex hulls and generalization
  bounds.
\newblock In \emph{International Conference on Computational Learning Theory},
  pages 59--73. Springer, 2002.

\bibitem[Cai and Zhang(2018)]{cai2018rate}
T~Tony Cai and Anru Zhang.
\newblock Rate-optimal perturbation bounds for singular subspaces with
  applications to high-dimensional statistics.
\newblock \emph{Ann. Statist.}, 46\penalty0 (1):\penalty0 60--89, 2018.

\bibitem[Cai et~al.(2013)Cai, Ma, and Wu]{cai2013sparse}
T~Tony Cai, Zongming Ma, and Yihong Wu.
\newblock Sparse pca: Optimal rates and adaptive estimation.
\newblock \emph{Ann. Statist.}, 41\penalty0 (6):\penalty0 3074--3110, 2013.

\bibitem[Cai et~al.(2015)Cai, Ma, and Wu]{Cai2015spiked}
T.~Tony Cai, Zongming Ma, and Yihong Wu.
\newblock Optimal estimation and rank detection for sparse spiked covariance
  matrices.
\newblock \emph{Probab. Theory Related Fields}, 161:\penalty0 781--815, 2015.

\bibitem[Cai et~al.(2016)Cai, Liang, and Rakhlin]{cai2016geometric}
T~Tony Cai, Tengyuan Liang, and Alexander Rakhlin.
\newblock Geometric inference for general high-dimensional linear inverse
  problems.
\newblock \emph{Ann. Statist.}, 44\penalty0 (4):\penalty0 1536--1563, 2016.

\bibitem[Candes and Plan(2011)]{candes2011tight}
Emmanuel~J Candes and Yaniv Plan.
\newblock Tight oracle inequalities for low-rank matrix recovery from a minimal
  number of noisy random measurements.
\newblock \emph{IEEE Trans. Inform. Theory}, 57\penalty0 (4):\penalty0
  2342--2359, 2011.

\bibitem[Chen and Cand{\`e}s(2018)]{chen2018projected}
Yuxin Chen and Emmanuel~J Cand{\`e}s.
\newblock The projected power method: An efficient algorithm for joint
  alignment from pairwise differences.
\newblock \emph{Comm. Pure Appl. Math.}, 71\penalty0 (8):\penalty0 1648--1714,
  2018.

\bibitem[Choi et~al.(2017)Choi, Taylor, and Tibshirani]{choi2017selecting}
Yunjin Choi, Jonathan Taylor, and Robert Tibshirani.
\newblock Selecting the number of principal components: Estimation of the true
  rank of a noisy matrix.
\newblock \emph{Ann. Statist.}, 45\penalty0 (6):\penalty0 2590--2617, 2017.

\bibitem[d'Aspremont et~al.(2005)d'Aspremont, Ghaoui, Jordan, and
  Lanckriet]{d2005direct}
Alexandre d'Aspremont, Laurent~E Ghaoui, Michael~I Jordan, and Gert~R
  Lanckriet.
\newblock A direct formulation for sparse pca using semidefinite programming.
\newblock In \emph{NIPS}, pages 41--48, 2005.

\bibitem[Deshpande and Montanari(2014)]{deshpande2014information}
Yash Deshpande and Andrea Montanari.
\newblock Information-theoretically optimal sparse pca.
\newblock In \emph{2014 IEEE Int. Symp. Info.}, pages 2197--2201. IEEE, 2014.

\bibitem[Deshpande et~al.(2014)Deshpande, Montanari, and
  Richard]{deshpande2014cone}
Yash Deshpande, Andrea Montanari, and Emile Richard.
\newblock Cone-constrained principal component analysis.
\newblock In \emph{NIPS}, pages 2717--2725, 2014.

\bibitem[Donoho and Gavish(2014)]{donoho2014minimax}
David Donoho and Matan Gavish.
\newblock Minimax risk of matrix denoising by singular value thresholding.
\newblock \emph{Ann. Statist.}, 42\penalty0 (6):\penalty0 2413--2440, 2014.

\bibitem[Donoho et~al.(2018)Donoho, Gavish, and Johnstone]{donoho2018optimal}
David~L Donoho, Matan Gavish, and Iain~M Johnstone.
\newblock Optimal shrinkage of eigenvalues in the spiked covariance model.
\newblock \emph{Ann. Statist.}, 46\penalty0 (4):\penalty0 1742, 2018.

\bibitem[d’Aspremont et~al.(2008)d’Aspremont, Bach, and
  Ghaoui]{d2008optimal}
Alexandre d’Aspremont, Francis Bach, and Laurent~El Ghaoui.
\newblock Optimal solutions for sparse principal component analysis.
\newblock \emph{J. Mach. Learn. Res.}, 9\penalty0 (Jul):\penalty0 1269--1294,
  2008.

\bibitem[Ferreira et~al.(2013)Ferreira, Iusem, and
  N{\'e}meth]{ferreira2013projections}
Orizon~Pereira Ferreira, Alfredo~N Iusem, and Sandor~Z N{\'e}meth.
\newblock Projections onto convex sets on the sphere.
\newblock \emph{Journal of Global Optimization}, 57\penalty0 (3):\penalty0
  663--676, 2013.

\bibitem[Giraud and Verzelen(2018)]{giraud2018partial}
Christophe Giraud and Nicolas Verzelen.
\newblock Partial recovery bounds for clustering with the relaxed $ k $ means.
\newblock \emph{arXiv preprint arXiv:1807.07547}, 2018.

\bibitem[Golub and Van~Loan(2012)]{golub2012matrix}
Gene~H Golub and Charles~F Van~Loan.
\newblock \emph{Matrix Computations}, volume~3.
\newblock JHU Press, 2012.

\bibitem[Haussler and Opper(1997{\natexlab{a}})]{haussler1997metric}
David Haussler and Manfred Opper.
\newblock Metric entropy and minimax risk in classification.
\newblock In \emph{Structures in Logic and Computer Science}, pages 212--235.
  Springer, 1997{\natexlab{a}}.

\bibitem[Haussler and Opper(1997{\natexlab{b}})]{haussler1997mutual}
David Haussler and Manfred Opper.
\newblock Mutual information, metric entropy and cumulative relative entropy
  risk.
\newblock \emph{Ann. Statist.}, 25\penalty0 (6):\penalty0 2451--2492,
  1997{\natexlab{b}}.

\bibitem[Javanmard et~al.(2016)Javanmard, Montanari, and
  Ricci-Tersenghi]{javanmard2016phase}
Adel Javanmard, Andrea Montanari, and Federico Ricci-Tersenghi.
\newblock Phase transitions in semidefinite relaxations.
\newblock \emph{P. Natl. Acad. Sci.}, 113\penalty0 (16):\penalty0 E2218--E2223,
  2016.

\bibitem[Jin and Wang(2016)]{jin2016influential}
Jiashun Jin and Wanjie Wang.
\newblock Influential features pca for high dimensional clustering.
\newblock \emph{Ann. Statist.}, 44\penalty0 (6):\penalty0 2323--2359, 2016.

\bibitem[Jin et~al.(2017)Jin, Ke, and Wang]{jin2017phase}
Jiashun Jin, Zheng~Tracy Ke, and Wanjie Wang.
\newblock Phase transitions for high dimensional clustering and related
  problems.
\newblock \emph{Ann. Statist.}, 45\penalty0 (5):\penalty0 2151--2189, 2017.

\bibitem[Johnstone(2001)]{johnstone2001distribution}
Iain~M Johnstone.
\newblock On the distribution of the largest eigenvalue in principal components
  analysis.
\newblock \emph{Ann. Statist.}, 29\penalty0 (2):\penalty0 295--327, 2001.

\bibitem[Journ{\'e}e et~al.(2010)Journ{\'e}e, Nesterov, Richt{\'a}rik, and
  Sepulchre]{journee2010generalized}
Michel Journ{\'e}e, Yurii Nesterov, Peter Richt{\'a}rik, and Rodolphe
  Sepulchre.
\newblock Generalized power method for sparse principal component analysis.
\newblock \emph{J. Mach. Learn. Res.}, 11\penalty0 (Feb):\penalty0 517--553,
  2010.

\bibitem[Kawale and Boley(2013)]{kawale2013constrained}
Jaya Kawale and Daniel Boley.
\newblock Constrained spectral clustering using l1 regularization.
\newblock In \emph{Proceedings of the 2013 SIAM International Conference on
  Data Mining}, pages 103--111. SIAM, 2013.

\bibitem[Kleindessner et~al.(2019)Kleindessner, Samadi, Awasthi, and
  Morgenstern]{kleindessner2019guarantees}
Matth{\"a}us Kleindessner, Samira Samadi, Pranjal Awasthi, and Jamie
  Morgenstern.
\newblock Guarantees for spectral clustering with fairness constraints.
\newblock \emph{arXiv preprint arXiv:1901.08668}, 2019.

\bibitem[Koltchinskii(2006)]{koltchinskii2006local}
Vladimir Koltchinskii.
\newblock Local rademacher complexities and oracle inequalities in risk
  minimization.
\newblock \emph{Ann. Statist.}, 34\penalty0 (6):\penalty0 2593--2656, 2006.

\bibitem[Krahmer et~al.(2014)Krahmer, Mendelson, and
  Rauhut]{krahmer2014suprema}
Felix Krahmer, Shahar Mendelson, and Holger Rauhut.
\newblock Suprema of chaos processes and the restricted isometry property.
\newblock \emph{Comm. Pure Appl. Math.}, 67\penalty0 (11):\penalty0 1877--1904,
  2014.

\bibitem[Lecu{\'e} and Mendelson(2009)]{lecue2009aggregation}
Guillaume Lecu{\'e} and Shahar Mendelson.
\newblock Aggregation via empirical risk minimization.
\newblock \emph{Probab. Theory Related Fields}, 145\penalty0 (3-4):\penalty0
  591--613, 2009.

\bibitem[L{\"o}ffler et~al.(2019)L{\"o}ffler, Zhang, and
  Zhou]{loffler2019optimality}
Matthias L{\"o}ffler, Anderson~Y Zhang, and Harrison~H Zhou.
\newblock Optimality of spectral clustering for gaussian mixture model.
\newblock \emph{arXiv preprint arXiv:1911.00538}, 2019.

\bibitem[Lu and Zhou(2016)]{lu2016statistical}
Yu~Lu and Harrison~H Zhou.
\newblock Statistical and computational guarantees of lloyd's algorithm and its
  variants.
\newblock \emph{arXiv preprint arXiv:1612.02099}, 2016.

\bibitem[Lugosi and Nobel(1999)]{lugosi1999adaptive}
G{\'a}bor Lugosi and Andrew~B Nobel.
\newblock Adaptive model selection using empirical complexities.
\newblock \emph{Ann. Statist.}, 27\penalty0 (6):\penalty0 1830--1864, 1999.

\bibitem[Ma et~al.(2019)Ma, Cai, and Li]{ma2019optimalb}
Rong Ma, T~Tony Cai, and Hongzhe Li.
\newblock Optimal and adaptive estimation of extreme values in the permuted
  monotone matrix model.
\newblock \emph{arXiv preprint arXiv:1911.12516}, 2019.

\bibitem[Ma et~al.(2020)Ma, Cai, and Li]{ma2019optimala}
Rong Ma, T~Tony Cai, and Hongzhe Li.
\newblock Optimal permutation recovery in permuted monotone matrix model.
\newblock \emph{J. Amer. Statist. Assoc.}, 2020.

\bibitem[Ma(2013)]{ma2013sparse}
Zongming Ma.
\newblock Sparse principal component analysis and iterative thresholding.
\newblock \emph{Ann. Statist.}, 41\penalty0 (2):\penalty0 772--801, 2013.

\bibitem[Massart(2007)]{massart2007concentration}
Pascal Massart.
\newblock \emph{Concentration Inequalities and Model Selection: Ecole d'Et{\'e}
  de Probabilit{\'e}s de Saint-Flour XXXIII-2003}.
\newblock Springer, 2007.

\bibitem[Montanari and Richard(2015)]{montanari2015non}
Andrea Montanari and Emile Richard.
\newblock Non-negative principal component analysis: Message passing algorithms
  and sharp asymptotics.
\newblock \emph{IEEE Trans. Inform. Theory}, 62\penalty0 (3):\penalty0
  1458--1484, 2015.

\bibitem[Ndaoud(2018)]{ndaoud2018sharp}
Mohamed Ndaoud.
\newblock Sharp optimal recovery in the two component gaussian mixture model.
\newblock \emph{arXiv preprint arXiv:1812.08078}, 2018.

\bibitem[Onaran and Villar(2017)]{onaran2017projected}
Efe Onaran and Soledad Villar.
\newblock Projected power iteration for network alignment.
\newblock In \emph{Wavelets and Sparsity XVII}, volume 10394, page 103941C.
  International Society for Optics and Photonics, 2017.

\bibitem[Paul(2007)]{paul2007asymptotics}
Debashis Paul.
\newblock Asymptotics of sample eigenstructure for a large dimensional spiked
  covariance model.
\newblock \emph{Statist. Sin.}, pages 1617--1642, 2007.

\bibitem[Perry et~al.(2018)Perry, Wein, Bandeira, and
  Moitra]{perry2018optimality}
Amelia Perry, Alexander~S Wein, Afonso~S Bandeira, and Ankur Moitra.
\newblock Optimality and sub-optimality of pca i: Spiked random matrix models.
\newblock \emph{Ann. Statist.}, 46\penalty0 (5):\penalty0 2416--2451, 2018.

\bibitem[Rakhlin et~al.(2017)Rakhlin, Sridharan, and
  Tsybakov]{rakhlin2017empirical}
Alexander Rakhlin, Karthik Sridharan, and Alexandre~B Tsybakov.
\newblock Empirical entropy, minimax regret and minimax risk.
\newblock \emph{Bernoulli}, 23\penalty0 (2):\penalty0 789--824, 2017.

\bibitem[Rangan and Fletcher(2012)]{rangan2012iterative}
Sundeep Rangan and Alyson~K Fletcher.
\newblock Iterative estimation of constrained rank-one matrices in noise.
\newblock In \emph{2012 IEEE Int. Symp. Info. Proceedings}, pages 1246--1250.
  IEEE, 2012.

\bibitem[Raskutti et~al.(2011)Raskutti, Wainwright, and
  Yu]{raskutti2011minimax}
Garvesh Raskutti, Martin~J Wainwright, and Bin Yu.
\newblock Minimax rates of estimation for high-dimensional linear regression
  over $\ell_q$-balls.
\newblock \emph{IEEE Trans. Inform. Theory}, 57\penalty0 (10):\penalty0
  6976--6994, 2011.

\bibitem[Shen and Huang(2008)]{shen2008sparse}
Haipeng Shen and Jianhua~Z Huang.
\newblock Sparse principal component analysis via regularized low rank matrix
  approximation.
\newblock \emph{J. Multiv. Anal.}, 99\penalty0 (6):\penalty0 1015--1034, 2008.

\bibitem[Singer(2011)]{singer2011angular}
Amit Singer.
\newblock Angular synchronization by eigenvectors and semidefinite programming.
\newblock \emph{Applied and Computational Harmonic Analysis}, 30\penalty0
  (1):\penalty0 20--36, 2011.

\bibitem[Szarek(1998)]{szarek1998metric}
Stanis{\l}aw Szarek.
\newblock Metric entropy of homogeneous spaces.
\newblock \emph{Banach Center Publications}, 43\penalty0 (1):\penalty0
  395--410, 1998.

\bibitem[Talagrand(2014)]{talagrand2014upper}
Michel Talagrand.
\newblock \emph{Upper and lower bounds for stochastic processes: modern methods
  and classical problems}, volume~60.
\newblock Springer Science \& Business Media, 2014.

\bibitem[Tsybakov(2009)]{tsybakov2009introduction}
Alexandre~B Tsybakov.
\newblock \emph{Introduction to Nonparametric Estimation}.
\newblock Springer Series in Statistics. Springer, New York, 2009.

\bibitem[Vershynin(2010)]{vershynin2010introduction}
Roman Vershynin.
\newblock Introduction to the non-asymptotic analysis of random matrices.
\newblock \emph{arXiv preprint arXiv:1011.3027}, 2010.

\bibitem[Vershynin(2018)]{vershynin2018high}
Roman Vershynin.
\newblock \emph{High-dimensional Probability: An Introduction with Applications
  in Data Science}, volume~47.
\newblock Cambridge University Press, 2018.

\bibitem[Verzelen(2012)]{verzelen2012minimax}
Nicolas Verzelen.
\newblock Minimax risks for sparse regressions: Ultra-high dimensional
  phenomenons.
\newblock \emph{Electronic Journal of Statistics}, 6:\penalty0 38--90, 2012.

\bibitem[Vu and Lei(2012)]{vu2012minimax}
Vincent Vu and Jing Lei.
\newblock Minimax rates of estimation for sparse pca in high dimensions.
\newblock In \emph{Artificial Intelligence and Statistics}, pages 1278--1286,
  2012.

\bibitem[Vu and Lei(2013)]{vu2013minimax}
Vincent~Q Vu and Jing Lei.
\newblock Minimax sparse principal subspace estimation in high dimensions.
\newblock \emph{Ann. Statist.}, 41\penalty0 (6):\penalty0 2905--2947, 2013.

\bibitem[Vu et~al.(2013)Vu, Cho, Lei, and Rohe]{vu2013fantope}
Vincent~Q Vu, Juhee Cho, Jing Lei, and Karl Rohe.
\newblock Fantope projection and selection: A near-optimal convex relaxation of
  sparse pca.
\newblock In \emph{NIPS}, pages 2670--2678, 2013.

\bibitem[Wang and Fan(2017)]{wang2017asymptotics}
Weichen Wang and Jianqing Fan.
\newblock Asymptotics of empirical eigenstructure for high dimensional spiked
  covariance.
\newblock \emph{Ann. Statist.}, 45\penalty0 (3):\penalty0 1342, 2017.

\bibitem[Wang and Davidson(2010)]{wang2010flexible}
Xiang Wang and Ian Davidson.
\newblock Flexible constrained spectral clustering.
\newblock In \emph{Proceedings of the 16th ACM SIGKDD International Conference
  on Knowledge Discovery and Data Mining}, pages 563--572. ACM, 2010.

\bibitem[Witten et~al.(2009)Witten, Tibshirani, and
  Hastie]{witten2009penalized}
Daniela~M Witten, Robert Tibshirani, and Trevor Hastie.
\newblock A penalized matrix decomposition, with applications to sparse
  principal components and canonical correlation analysis.
\newblock \emph{Biostatistics}, 10\penalty0 (3):\penalty0 515--534, 2009.

\bibitem[Wu and Yang(2016)]{wu2016minimax}
Yihong Wu and Pengkun Yang.
\newblock Minimax rates of entropy estimation on large alphabets via best
  polynomial approximation.
\newblock \emph{IEEE Trans. Inform. Theory}, 62\penalty0 (6):\penalty0
  3702--3720, 2016.

\bibitem[Yang et~al.(2011)Yang, Ma, and Buja]{yang2011sparse}
Dan Yang, Zongming Ma, and Andreas Buja.
\newblock A sparse svd method for high-dimensional data.
\newblock \emph{arXiv preprint arXiv:1112.2433}, 2011.

\bibitem[Yang(1999)]{yang1999minimax}
Yuhong Yang.
\newblock Minimax nonparametric classification. i. rates of convergence.
\newblock \emph{IEEE Trans. Inform. Theory}, 45\penalty0 (7):\penalty0
  2271--2284, 1999.

\bibitem[Yang and Barron(1999)]{yang1999information}
Yuhong Yang and Andrew Barron.
\newblock Information-theoretic determination of minimax rates of convergence.
\newblock \emph{Ann. Statist.}, pages 1564--1599, 1999.

\bibitem[Yatracos(1988)]{yatracos1988lower}
Yannis~G Yatracos.
\newblock A lower bound on the error in nonparametric regression type problems.
\newblock \emph{Ann. Statist.}, pages 1180--1187, 1988.

\bibitem[Yuan and Zhang(2013)]{yuan2013truncated}
Xiao-Tong Yuan and Tong Zhang.
\newblock Truncated power method for sparse eigenvalue problems.
\newblock \emph{J. Mach. Learn. Res.}, 14\penalty0 (Apr):\penalty0 899--925,
  2013.

\bibitem[Zhang et~al.(2018)Zhang, Cai, and Wu]{zhang2018heteroskedastic}
Anru Zhang, T~Tony Cai, and Yihong Wu.
\newblock Heteroskedastic pca: Algorithm, optimality, and applications.
\newblock \emph{arXiv preprint arXiv:1810.08316}, 2018.

\bibitem[Zou et~al.(2006)Zou, Hastie, and Tibshirani]{zou2006sparse}
Hui Zou, Trevor Hastie, and Robert Tibshirani.
\newblock Sparse principal component analysis.
\newblock \emph{J. Comput. Graph. Stat.}, 15\penalty0 (2):\penalty0 265--286,
  2006.

\end{thebibliography}

\end{document}